\documentclass[a4paper,12pt]{article}

\usepackage{amsmath,amssymb,latexsym,amsfonts,amsthm}
\usepackage{mathrsfs}

\theoremstyle{plain}
\newtheorem{Thm}{Theorem}[section]

\newtheorem{Lem}[Thm]{Lemma}
\newtheorem{Prop}[Thm]{Proposition}
\newtheorem{Cor}[Thm]{Corollary}

\theoremstyle{definition}

\newtheorem{Rem}[Thm]{Remark}
\newtheorem{Ex}[Thm]{Example}

\newcommand{\bD}{\ensuremath{\mathbb{D}}}
\newcommand{\bE}{\ensuremath{\mathbb{E}}}

\newcommand{\bR}{\ensuremath{\mathbb{R}}}

\newcommand{\cC}{\ensuremath{\mathcal{C}}}

\newcommand{\cE}{\ensuremath{\mathcal{E}}}
\newcommand{\cF}{\ensuremath{\mathcal{F}}}

\newcommand{\sM}{\ensuremath{\mathscr{M}}}

\newcommand{\sP}{\ensuremath{\mathscr{P}}}

\newcommand{\sW}{\ensuremath{\mathscr{W}}}

\newcommand{\vn}{\ensuremath{\mbox{{\boldmath $n$}}}}

\newcommand{\eps}{\ensuremath{\varepsilon}}
\newcommand{\e}{{\rm e}}
\renewcommand{\d}{{\rm d}}

\newcommand{\claw}{\stackrel{{\rm law}}{\longrightarrow}}

\renewcommand{\hat}{\widehat}
\renewcommand{\tilde}{\widetilde}
\renewcommand{\bar}{\overline}

\newcommand{\abra}[1]{\left| #1 \right|}
\newcommand{\cbra}[1]{\left( #1 \right)}
\newcommand{\kbra}[1]{\left\{ #1 \right\}}
\newcommand{\ebra}[1]{\left[ #1 \right]}

\numberwithin{equation}{section}

\newcounter{No}

\newcounter{Ci}[subsection]

\makeatletter
\renewcommand\section{\@startsection {section}{1}{\z@}%
                                   {-3.5ex \@plus -1ex \@minus -.2ex}%
                                   {2.3ex \@plus.2ex}%
                                   {\normalfont\large\bf}}
\makeatother

\makeatletter
\@addtoreset{footnote}{page}
\makeatother

\theoremstyle{plain}

\renewcommand{\dagger}{h}

\newcounter{CO}
\newcommand{\Confirst}[1]{\refstepcounter{CO} \label{con: #1}}

\newcounter{EO}

\usepackage{graphicx}
\newcommand{\mail}{
\scalebox{0.6}{\includegraphics{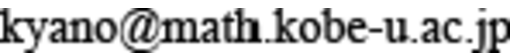}}
}

\begin{document}

%
\begin{center}
{\Large \bf 
Penalising symmetric stable L\'evy paths 
}
\end{center}
\begin{center}
Kouji \textsc{Yano}\footnote{
Department of Mathematics, Graduate School of Science, Kobe University, Kobe, Japan.\\
E-mail: \mail}\footnote{
The research of this author is supported by KAKENHI (20740060)}, \qquad 
Yuko \textsc{Yano}\footnote{
Research Institute for Mathematical Sciences, Kyoto University, Kyoto, Japan.\label{foot: RIMS}} 
\qquad and \qquad 
Marc \textsc{Yor}\footnote{
Laboratoire de Probabilit\'es et Mod\`eles Al\'eatoires, Universit\'e Paris VI, Paris, France.}\footnote{
Institut Universitaire de France}\footnotemark[3]
\end{center}
\begin{center}
{\small \today}
\end{center}
\bigskip


\begin{abstract}
Limit theorems for the normalized laws with respect to two kinds of weight functionals 
are studied for any symmetric stable L\'evy process of index $ 1 < \alpha \le 2 $. 
The first kind is a function of the local time at the origin, 
and the second kind is the exponential of an occupation time integral. 
Special emphasis is put on the role played by a stable L\'evy counterpart 
of the universal $ \sigma $-finite measure, 
found in \cite{NRY} and \cite{NRY2}, 
which unifies the corresponding limit theorems in the Brownian setup 
for which $ \alpha =2 $. 
\end{abstract}

\section{Introduction}

Roynette, Vallois and Yor (\cite{MR2261065}, 
\cite{MR2229621} 
and \cite{MR2253307} 
and references therein) 
have shown the existence of the limit laws for normalized Wiener measures 
with respect to various weight processes; 
we call these studies {\em penalisation problems}. 
Najnudel, Roynette and Yor (see \cite{RY}, \cite{MR2339327}, \cite{NRY} and \cite{NRY2}) 
have recently discovered that 
these penalisation problems may be unified 
with the help of the following ``universal" $ \sigma $-finite measure on the canonical space: 
\begin{align}
\sW = \int_0^{\infty } \frac{\d u}{\sqrt{2 \pi u}} W^{(u)} \bullet P^{\rm 3B}_0 
\label{eq: curly W}
\end{align}
where $ W^{(u)} $ stands for the law of the brownian bridge from $ 0 $ to $ 0 $ of length $ u $, 
$ P^{\rm 3B}_0 $ for that of the symmetrized 3-dimensional Bessel process starting from 0, 
i.e., $ P^{\rm 3B}_0 = (P^{\rm 3B,+}_0 + P^{\rm 3B,-}_0)/2 $, 
and the symbol $ \bullet $ for the concatenation between the laws of these two processes. 

The purpose of the present paper 
is to develop some of these penalisation problems 
in the case of any symmetric stable L\'evy process of index $ 1<\alpha \le 2 $. 
As an analogue of $ \sW $, 
we introduce the following $ \sigma $-finite measure 
\begin{align}
\sP = \frac{\Gamma(1/\alpha )}{\alpha \pi} 
\int_0^{\infty } \frac{\d u}{u^{1/\alpha }} Q^{(u)} \bullet P^{\dagger}_0 
\label{eq: sP intro}
\end{align}
where 
$ Q^{(u)} $ stands for the law of the bridge from $ 0 $ to $ 0 $ of length $ u $ 
and 
$ P^{\dagger}_0 $ for the $ h $-path process of the killed process 
with respect to the function $ |x|^{\alpha -1} $. 
We shall put some special emphasis on the role played 
by the universal $ \sigma $-finite measure $ \sP $ 
which helps to unify our penalisation problems.

Let $ \bD $ denote the canonical space of c\`adl\`ag paths $ w:[0,\infty ) \to \bR $. 
Let $ (X_t) $ denote the coordinate process, 
$ (\cF_t) $ its natural filtration, 
and $ \cF_{\infty } = \vee_{t \ge 0} \cF_t $. 
Let $ (P_x) $ denote the law on $ \bD $ 
of the symmetric stable process of index $ 1<\alpha \le 2 $ 
such that $ P_0[\e^{i \lambda X_t}] = \e^{- t |\lambda|^{\alpha }} $ for $ \lambda \in \bR $. 
Note that, if $ \alpha =2 $, then $ (X_t) $ has the same law 
as $ \sqrt{2} $ times the standard brownian motion. 

We say that {\em a family of measures $ \{ \sM_t \}_{t \ge 0} $ on $ \cF_{\infty } $ 
converges as $ t \to \infty $ to a measure $ \sM $ along $ (\cF_s) $} 
if, for each $ s>0 $, we have 
$ \sM_t[Z_s] \to \sM[Z_s] $ as $ t \to \infty $ 
for all bounded $ \cF_s $-measurable functionals $ Z_s $. 
For a measure $ \sM $ on $ \cF_{\infty } $ and a functional $ F $ 
measurable with respect to $ \cF_{\infty } $, 
the symbol $ F \cdot \sM $ stands for the measure $ A \mapsto \sM[1_A F] $. 
Let $ x \in \bR $ be fixed. 
Then penalisation problems are stated as follows: 

\begin{quote}
{\bf Question 1.} 
Let $ \Gamma = (\Gamma_t:t \ge 0) $ be a given non-negative process 
such that $ P_x[\Gamma_t] \neq 0 $ for large enough $ t $. 
\\ 
{\bf (Q1)} 
Does there exist a limit probability measure $ P^{\Gamma}_x $ 
such that 
\begin{align}
\frac{\Gamma_t \cdot P_x}{P_x[\Gamma_t]} 
\ \stackrel{t \to \infty }{\longrightarrow} \ 
P^{\Gamma}_x 
\qquad \text{along $ (\cF_s) $?} 
\label{eq: Q1}
\end{align}
{\bf (Q2)} 
How can one characterise the limit probability measure $ P^{\Gamma}_x $ assuming it exists? 
\end{quote}

For each $ x \in \bR $, let $ \sP_x $ denote the law of $ (x+X_t:t \ge 0) $ under $ \sP $. 
We can gain a clear insight into some of these penalisation problems 
if we answer the following 

\begin{quote}
{\bf Question 2.} 
Let $ \Gamma $ as above. 
\\ 
{\bf (Q1$ ' $)} 
Can one find a positive function $ \mu(t) $ 
and a measurable functional $ \Gamma_{\infty } $ 
such that 
\begin{align}
\frac{\Gamma_t \cdot P_x}{\mu(t)} 
\ \stackrel{t \to \infty }{\longrightarrow} \ 
\Gamma_{\infty } \cdot \sP_x 
\qquad \text{along $ (\cF_s) $?} 
\label{eq: Q1'}
\end{align}
{\bf (Q2$ ' $)} 
For any non-negative $ \sP_x $-integrable functional $ F $, 
can one find a non-negative $ (\cF_t,P_x) $-martingale $ (M_{t,x}(F):t \ge 0) $ such that 
\begin{align}
(F \cdot \sP_x) |_{\cF_t} = M_{t,x}(F) \cdot P_x|_{\cF_t} 
, \qquad t \ge 0 ? 
\label{eq: Q2'}
\end{align}
\end{quote}

If we can find such a function $ \mu(t) $ as in \eqref{eq: Q1'} 
and if $ 0 < \sP_x[\Gamma_{\infty }] < \infty $, 
then we obtain the convergence \eqref{eq: Q1} with the limit probability measure 
\begin{align}
P^{\Gamma}_x = \frac{\Gamma_{\infty } \cdot \sP_x}{\sP_x[\Gamma_{\infty }]} . 
\label{}
\end{align}
We shall prove in Theorem \ref{thm: mart op} 
that there exist such martingales $ (M_{t,x}(F)) $ as in \eqref{eq: Q2'}. 
We shall call $ M_{t,x}(\cdot) $ the {\em martingale generator} 
and we shall study its properties in Sections \ref{sec: univ} and \ref{sec: further}. 
Then the limit probability measure $ P^{\Gamma}_x $ is characterised by 
\begin{align}
P^{\Gamma}_x |_{\cF_t} = \frac{M_{t,x}(\Gamma_{\infty })}{\sP_x[\Gamma_{\infty }]} \cdot P_x|_{\cF_t} 
, \qquad t \ge 0 . 
\label{}
\end{align}
Therefore, if we answer Question 2, then we have answered Question 1. 

Our strategy to answer {\bf (Q1$ ' $)} is as follows. 
Since the index $ \alpha $ is supposed to be in $ (1,2] $, 
each point of $ \bR $ is regular and recurrent. 
Hence, associated with the process, 
there is a jointly continuous local time $ (L(t,x)) $. 
We simply write $ L_t = L(t,0) $, 
and, associated with this local time, 
there is {\em It\^o's measure $ \vn $ of excursions away from the origin} 
(see Section \ref{sec: exc}). 
Let $ R $ denote the {\em lifetime} of an excursion path. 
For $ t>0 $, we define $ M^{(t)} $ as the probability measure on $ \cF_t $ given by 
\begin{align}
M^{(t)} = \frac{1_{\{ R>t \}}}{\vn(R>t)} \cdot \vn |_{\cF_t} 
\label{eq: Meander}
\end{align}
and here we call $ M^{(t)} $ the distribution of the {\it stable meander}. 
We remark that 
our meander distribution \eqref{eq: Meander} 
is definitely different from that of \cite{MR1465814} etc. 
where the meander is defined by conditioning on $ \{ R>t \} $ 
the excursion process for the {\em reflected} stable L\'evy process 
$ (X_t-\min_{s \le t} X_s:t \ge 0) $. 
We shall prove the following formula (Theorem \ref{thm: LED}) 
of disintegration of $ P_0|_{\cF_t} $ for each $ t>0 $ 
with respect to last exit time from the origin: 
\begin{align}
\frac{1}{\vn(R>t)} P_0|_{\cF_t} 
= 
\frac{\Gamma(1/\alpha )}{\alpha \pi} 
\int_0^t \cbra{ 1-\frac{u}{t} }^{\frac{1}{\alpha }-1} 
\frac{\d u}{u^{1/\alpha }} Q^{(u)} \bullet M^{(t-u)} . 
\label{eq: LED2}
\end{align}
As a check, the total masses of both sides agree, 
as we shall show in Proposition \ref{prop: vn Rt}. 
Then, we shall establish (in Theorem \ref{thm: meander}) the convergence 
\begin{align}
M^{(t)} 
\stackrel{t \to \infty }{\longrightarrow} 
P^{\dagger}_0 
\qquad \text{along $ (\cF_s) $.} 
\label{eq: meander conv Pdagger2}
\end{align}
Noting that $ (1-\frac{u}{t})^{\frac{1}{\alpha } -1} \to 1 $ as $ t \to \infty $, 
we may expect that, in some sense: 
\begin{align}
\int_0^t \cbra{ 1-\frac{u}{t} }^{\frac{1}{\alpha }-1} 
\frac{\d u}{u^{1/\alpha }} Q^{(u)} \bullet M^{(t-u)} 
\stackrel{t \to \infty }{\longrightarrow} 
\int_0^{\infty } \frac{\d u}{u^{1/\alpha }} Q^{(u)} \bullet P^{\dagger}_0 . 
\label{eq: convergence of P_0/nRt to sP}
\end{align}
We shall prove several analytic lemmas 
which justify the convergence \eqref{eq: convergence of P_0/nRt to sP} 
and then 
we shall establish the convergence \eqref{eq: Q1'} 
with the function $ \mu(t) = \vn(R>t) $. 

In order to answer Question 2 (and in particular {\bf (Q2$ ' $)}), 
we shall establish the convergence \eqref{eq: Q1'} 
and compute the martingale generator 
by case study. 
We confine ourselves to the following two kinds of weight functionals: 
\begin{quote}
(i) 
$ \Gamma_t = f(L_t) $ for some non-negative Borel functions $ f $ 
with some integrability property; 
\\
(ii) 
$ \displaystyle \Gamma_t = \exp \kbra{ - \int L(t,x) V(\d x) } $ 
for some non-negative Borel measure $ V $. 
We call the problems in such a case the {\em Feynman--Kac penalisations}. 
\end{quote}

The organisation of the present paper is as follows. 
In Section \ref{sec: prel} we recall some preliminary facts 
about symmetric stable L\'evy processes. 
In Section \ref{sec: exc} we study It\^o's measure of excursions away from the origin 
relatively to the symmetric stable process. 
In Section \ref{sec: intro2} 
we prove several formulae concerning the stable meander and $ h $-path process, 
which play important roles in the study of our penalisation problems. 
In Section \ref{sec: univ} we make 
general observations on the universal $ \sigma $-finite measure $ \sP_x $ 
and the martingale generator $ M_{t,x}(\cdot) $. 
In Section \ref{sec: imp lem} we prove several convergence lemmas 
which play fundamental roles in the proof of our penalisation problems. 
Section \ref{sec: LT penal} is devoted to the study of penalisations 
with a function of the local time at the origin. 
Section \ref{sec: FK penal} is devoted to the study of Feynman--Kac penalisations. 
In Section \ref{sec: further} we characterise certain non-negative 
$ (P_0,\cF_t) $-martingales in terms of $ \sP $.


\section{Preliminaries about the symmetric stable process of index $ 1 < \alpha \le 2 $}
\label{sec: prel}

Recall that $ (X_t,\cF_t,P_x) $ is the canonical representation 
of a one-dimensional symmetric stable L\'evy process of index $ 1<\alpha \le 2 $ 
such that 
\begin{align}
P_0[\e^{i \lambda X_t}] = \e^{-t |\lambda|^{\alpha }} 
\quad \text{for} \ 
\lambda \in \bR . 
\label{}
\end{align}
All results presented in this section are well-known; 
see, e.g., \cite{MR1406564}. 

\noindent 
{\bf 1).} 
$ (X_t) $ has a transition density 
$ P_x(X_t \in \d y) = p_t(y-x) \d y $ 
where 
$ p_t(x) $ is given by 
\begin{align}
p_t(x) 
= \frac{1}{\pi} \int_0^{\infty } (\cos x \lambda) \e^{-t \lambda^{\alpha }} \d \lambda 
. 
\label{}
\end{align}
For $ q>0 $, we set 
\begin{align}
u_q(x) = \int_0^{\infty } \e^{-qt} p_t(x) \d t 
= \frac{1}{\pi} \int_0^{\infty } \frac{\cos x \lambda}{q+\lambda^{\alpha }} \d \lambda . 
\label{eq: uqx}
\end{align}
In particular, if we take $ x=0 $, we have 
\begin{align}
p_t(0) = p_1(0) t^{-\frac{1}{\alpha }} 
\quad \text{where} \ 
p_1(0) = \frac{\Gamma (1/\alpha )}{\alpha \pi} 
\label{}
\end{align}
and 
\begin{align}
u_q(0) = u_1(0) q^{\frac{1}{\alpha }-1} 
\quad \text{where} \ 
u_1(0) = \frac{\Gamma(1-1/\alpha ) \Gamma(1/\alpha )}{\alpha \pi} . 
\label{eq: uq0}
\end{align}

\noindent 
{\bf 2).} 
Let $ T_{\{ a \}} $ denote the first hitting time of $ a $ for the coordinate process $ (X_t) $: 
\begin{align}
T_{\{ a \}} = \inf \{ t>0: X_t=a \} . 
\label{}
\end{align}
Then the Laplace transform of the law of $ T_{\{ 0 \}} $ is given by 
\begin{align}
P_x[\e^{-q T_{\{ 0 \}}}] = \frac{u_q(x)}{u_q(0)} 
, \qquad x \in \bR , \ q>0 
\label{eq: LT of T0}
\end{align}
(see, e.g., \cite[pp. 64]{MR1406564}). 
For further study of the law of $ T_{\{ 0 \}} $, see \cite{YYY2}. 

Since $ T_{\{ y \}} $ under $ P_x $ has the same law as $ T_{\{ 0 \}} $ under $ P_{x-y} $, 
the formula \eqref{eq: LT of T0} implies the following facts: 
\\ \quad 
{\rm (i)} 
Each point is a recurrent state, i.e., 
$ P_x(T_{\{ y \}}<\infty ) = 1 $ for any $ x,y \in \bR $ with $ x \neq y $; 
\\ \quad 
{\rm (ii)} 
Each point is regular for itself, i.e., 
$ P_x(T_{\{ x \}}=0)=1 $ for any $ x \in \bR $.

\noindent 
{\bf 3).} 
The process admits a jointly continuous local time $ L(t,x) $ such that 
\begin{align}
L(t,x) = \lim_{\eps \to 0+} \frac{1}{2 \eps} \int_0^t 1_{\{ |X_s-x|<\eps \}} \d s 
\label{}
\end{align}
almost surely. 
We simply write $ L_t = L(t,0) $. 
Denote the inverse local time at the origin by $ \tau_l = \inf \{ t>0: L_t>l \} $. 
Then $ (\tau_l:l \ge 0) $ is a stable subordinator of index $ 1-1/\alpha $ 
such that 
\begin{align}
P_0[\e^{-q \tau_l}] = \e^{-l/u_q(0)} 
\label{}
\end{align}
(see, e.g., \cite[pp. 131]{MR1406564}), 
where $ u_q(0) $ is given explicitly by \eqref{eq: uq0}. 
Let $ \theta_t:\bD \to \bD $ stand for the shift operator: 
$ \theta_t(w) = w(t+\cdot) $. 
Since $ \tau_l = T_{\{ 0 \}} + \tau_l \circ \theta_{T_{\{ 0 \}}} $, we have 
\begin{align}
P_x \ebra{ \int_0^{\infty } \e^{-qt} \d L_t } 
=& P_x \ebra{ \int_0^{\infty } \e^{-q \tau_l} \d l } 
= P_x[\e^{-q T_{\{ 0 \}}}] \int_0^{\infty } P_0[\e^{-q\tau_l}] \d l 
\\
=& \frac{u_q(x)}{u_q(0)} \cdot u_q(0) 
= \int_0^{\infty } \e^{-qt} p_t(x) \d t 
\label{}
\end{align}
for all $ q>0 $. 
Hence we see that 
\begin{align}
P_x \ebra{ \int_0^{\infty } f(t) \d L_t } 
= \int_0^{\infty } f(t) p_t(x) \d t 
\label{}
\end{align}
for any non-negative measurable function $ f $ on $ [0,\infty ) $. 
Consequently, we may write 
\begin{align}
P_x[\d L_t] = p_t(x) \d t 
, \qquad x \in \bR. 
\label{}
\end{align}

\section{It\^o's measure of excursions away from the origin}
\label{sec: exc}

Since the origin is a regular and recurrent state, 
we can apply It\^o's excursion theory 
(\cite{MR0402949}; see also \cite{MR1406564} and \cite{MR1138461} for details). 

We denote by $ \bE $ the set of c\`adl\`ag paths $ e:[0,\infty ) \to \bR \cup \{ \Delta \} $ 
such that 
\begin{align}
\begin{cases}
e(t) \in \bR \setminus \{ 0 \} 
\quad & \text{for $ 0<t<R(e) $}, \\
e(t)=\Delta 
\quad & \text{for $ t \ge R(e) $} 
\end{cases}
\label{}
\end{align}
where 
\begin{align}
R = R(e) = \inf \{ t>0: e(t)=\Delta \} . 
\label{}
\end{align}
We call $ \bE $ {\em the set of excursions} 
and every element $ e $ of $ \bE $ an {\em excursion path}. 
For an excursion path $ e \in \bE $, 
we call $ R(e) $ the {\em lifetime} of $ e $. 
The point $ \Delta $ is called the {\em cemetery}. 

We set $ D = \{ l: \tau_l - \tau_{l-}>0 \} $. 
For each $ l \in D $, 
we set 
\begin{align}
e_l(t) = 
\begin{cases}
X_{t+\tau_{l-}} 
, \quad & \text{for} \ 0 \le t < \tau_l - \tau_{l-} , \\
\Delta 
, \quad & \text{for} \ t \ge \tau_l - \tau_{l-}. 
\end{cases}
\label{}
\end{align}
Then It\^o's fundamental theorem (\cite{MR0402949}) asserts that 
the point process $ (e_l: l \in D) $ taking values on $ \bE $ 
is a Poisson point process. 
Its characteristic measure will be denoted by $ \vn $ 
and called {\em It\^o's measure of excursions away from the origin}. 
It\^o's measure $ \vn $ is a $ \sigma $-finite measure on any $ \cF_t $ 
which has no mass outside the set 
\begin{align}
\{ e \in \bE : X_0(e)=0 , \ 0<R(e)<\infty \} . 
\label{}
\end{align}
For the fact that $ \vn(\{ X_0 = 0 \}^c) = 0 $, see \cite{Y}.

For $ x \in \bR \setminus \{ 0 \} $, 
we denote by $ P^0_x $ 
the law of the killed process, i.e., 
the law on $ \bE $ of the path $ (X^0_t) $ under $ P_x $ where 
\begin{align}
X^0_t = 
\begin{cases}
X_t 
, \quad & 0 \le t < T_{\{ 0 \}} , \\
\Delta 
, \quad & t \ge T_{\{ 0 \}} . 
\end{cases}
\label{}
\end{align}
We shall utilise the following formulae.

\begin{Thm}[Markov property of $ \vn $] \label{thm: Markov of vn}
It holds that 
\begin{align}
\vn[Z_t F(X \circ \theta_t)] 
= \int \vn[Z_t;X_t \in \d x] P^0_x[F(X)] 
\label{}
\end{align}
for any $ t>0 $, 
any non-negative $ \cF_t $-measurable functional $ Z_t $ 
and any non-negative measurable functional $ F $ on $ \bE $. 
\end{Thm}

\begin{Thm}[Compensation formula] 
\label{thm: compensation}
Let $ F=F(t,\omega,e) $ be a measurable functional on 
$ [0,\infty ) \times \bD \times \bE $ such that, 
for every fixed $ e \in \bE $, 
the process $ (F(t,\cdot,e):t \ge 0) $ is $ (\cF_t) $-predictable. 
Then\footnote{Here the symbol $ \tilde{} $ means independence.} 
\begin{align}
P_0 \ebra{ \sum_{l \in D} F(\tau_{l-},X,e_l) } 
= P_0 \otimes \tilde{\vn} \ebra{ \int_0^{\infty } \d L_t F(t,X,\tilde{X}) } . 
\label{eq: compensation}
\end{align}
\end{Thm}

We omit the proofs of Theorems \ref{thm: Markov of vn} and \ref{thm: compensation}. 
For their proofs, see \cite{MR1406564}, \cite{MR1138461} and \cite{MR1725357}.

\subsection{Entrance law}

In order to characterise the entrance law, we need the following 

\begin{Thm}[\cite{CFY} and \cite{MR2247835}] \label{thm: master1}
For any non-negative measurable function $ f $ on $ \bR $, it holds that 
\begin{align}
\int_0^{\infty } \e^{-qt} \vn \ebra{f(X_t)} \d t 
= \int f(x) P_x \ebra{ \e^{-q T_{\{ 0 \}}} } \d x . 
\label{eq: master1}
\end{align}
\end{Thm}

We remark that the relation \eqref{eq: master1} can be found 
in Chen--Fukushima--Ying \cite[Eq. (2.8)]{CFY} 
and Fitzsimmons--Getoor \cite[Eq. (3.22)]{MR2247835} 
in a fairly general Markovian framework as 
\begin{align}
\int_0^{\infty } \e^{-qt} \vn[f(X_t)] \d t 
= \int f(x) \hat{P}_x \ebra{ \e^{-q T_{\{ 0 \}}(\hat{X})} } m(\d x) 
\label{eq: master CFY and FG}
\end{align}
where 
$ (X_t,P_x) $ and $ (\hat{X},\hat{P}_x) $ are 
in weak duality with respect to the reference measure $ m $. 
In our case, $ (\hat{X},\hat{P}_x) = (-X_t,P_x) $ 
and $ m(\d x) = \d x $, the Lebesgue measure. 
Although \eqref{eq: master1} is a special case of \eqref{eq: master CFY and FG}, 
we give the proof of Theorem \ref{thm: master1} for completeness of this paper. 

\begin{proof}[Proof of Theorem \ref{thm: master1}]
Note that 
\begin{align}
\int_0^{\infty } \e^{-qt} f(X_t) \d t 
= \sum_{l \in D} \e^{-q \tau_{l-}} \int_0^{R(e_l)} \e^{-qt} f(e_l(t)) \d t . 
\label{}
\end{align}
By Theorem \ref{thm: compensation}, we obtain 
\begin{align}
P_0 \ebra{ \int_0^{\infty } \e^{-qt} f(X_t) \d t } 
= P_0 \ebra{ \int_0^{\infty } \e^{-qt} \d L_t } \vn \ebra{ \int_0^R \e^{-qt} f(X_t) \d t } . 
\label{}
\end{align}
Since $ P_0 \ebra{ \int_0^{\infty } \e^{-qt} \d L_t } = u_q(0) $, 
we have 
\begin{align}
\int_0^{\infty } \e^{-qt} \vn \ebra{f(X_t)} \d t 
= \int f(x) \frac{u_q(x)}{u_q(0)} \d x . 
\label{eq: master1-1}
\end{align}
By the identity \eqref{eq: LT of T0}, 
we obtain \eqref{eq: master1}. 
The proof is complete. 
\end{proof}

The following formula holds: 

\begin{Prop} \label{prop: vn Rt}
\begin{align}
\vn(R>t) = \vn(R>1) t^{\frac{1}{\alpha }-1} 
\label{}
\end{align}
where 
\begin{align}
\vn(R>1) 
= 
\frac{\alpha \pi}{\Gamma(1-1/\alpha) \Gamma (1/\alpha )^2} . 
\label{}
\end{align}
In particular, 
\begin{align}
\frac{\vn(R>t-s)}{\vn(R>t)} = \cbra{ 1-\frac{s}{t} }^{\frac{1}{\alpha }-1} 
\qquad \text{for} \ 
0<s<t . 
\label{eq: con R ratio}
\end{align}
\end{Prop}

Although it is well-known, we again give the proof for completeness of this paper. 

\begin{proof}
Take $ f = 1 $ in \eqref{eq: master1-1}. 
Then we have 
$ \vn[f(X_t)] = \vn(R>t) $, 
and the identity \eqref{eq: master1-1} implies that 
\begin{align}
\int_0^{\infty } \e^{-qt} \vn(R>t) \d t 
= \frac{1}{q u_q(0)} = \frac{1}{u_1(0)} q^{-1/\alpha } . 
\label{eq: LT vnR}
\end{align}
This completes the proof. 
\end{proof}

The following theorem characterises the entrance law. 

\begin{Thm} \label{thm: rho}
There exists a bi-measurable function $ \rho(t,x) $ 
which is at the same time 
a space density of the entrance law 
\begin{align}
\vn(X_t \in \d x) = \rho(t,x) \d x 
\label{eq: ent law density}
\end{align}
and a time density of the first hitting time 
\begin{align}
P_x(T_{\{ 0 \}} \in \d t) = \rho(t,x) \d t . 
\label{eq: first pas density}
\end{align}
That is, 
\begin{align}
\rho(t,x) = \frac{\vn(X_t \in \d x)}{\d x} = \frac{P_x(T_{\{ 0 \}} \in \d t)}{\d t} . 
\label{eq: ent law and first pas}
\end{align}
\end{Thm}

\begin{proof}
Note that 
$ P^0_x(X_t \in \d y) = p^0_t(x,y) \d y $ where 
\begin{align}
p^0_t(x,y) = p_t(y-x) - \int_0^t p_{t-s}(y) P_x(T_{\{ 0 \}} \in \d s) . 
\label{}
\end{align}
Now we set 
\begin{align}
\rho(t,x) = \int \vn(X_{t/2} \in \d y) p^0_{t/2}(y,x) . 
\label{}
\end{align}
Let $ f $ be a non-negative measurable function on $ \bR $. 
By the Markov property, 
we see that 
\begin{align}
\vn[f(X_t)] = \int \vn(X_{t/2} \in \d y) P^0_y[f(X_{t/2})] 
= \int f(x) \rho(t,x) \d x . 
\label{eq: ent law integral}
\end{align}
Hence we obtain \eqref{eq: ent law density}. 
Using the formulae \eqref{eq: ent law integral} 
and \eqref{eq: master1}, 
we see that 
\begin{align}
\int \d x f(x) \int_0^{\infty } \e^{-qt} \rho(t,x) \d t 
=& 
\int_0^{\infty } \e^{-qt} \vn[f(X_t)] \d t 
\\
=& 
\int \d x f(x) P_x[\e^{-q T_{\{ 0 \}}}] . 
\label{}
\end{align}
Hence we obtain \eqref{eq: first pas density}. 
\end{proof}

\section{Stable meander and $ h $-path process} \label{sec: intro2}

\subsection{Disintegration with respect to the last exit time}

For $ u>0 $, let $ Q^{(u)} $ denote the law of the bridge 
$ P_0(\cdot | X_u=0) $ 
considered to be a probability measure on $ \cF_u $. 
We denote by $ X^{(u)} = (X_t:0 \le t \le u) $ 
the coordinate process considered up to time $ u $. 
We denote the {\em concatenation} between 
the two processes $ X^{(u)} $ and $ \tilde{X}^{(v)}=(\tilde{X}_t:0 \le t \le v) $ 
by $ X^{(u)} \bullet \tilde{X}^{(v)} = ((X^{(u)} \bullet \tilde{X}^{(v)})_t:0 \le t \le u+v) $: 
\begin{align}
\cbra{ X^{(u)} \bullet \tilde{X}^{(v)} }_t = 
\begin{cases}
X^{(u)}_t 
, \quad & 0 \le t < u , \\
\tilde{X}^{(v)}_{t-u} 
, \quad & u \le t \le u+v . 
\end{cases}
\label{}
\end{align}
The measure $ Q^{(u)} \bullet M^{(v)} $ 
is defined as the law of the concatenation $ X^{(u)} \bullet \tilde{X}^{(v)} $ 
between the two processes $ X^{(u)} $ and $ \tilde{X}^{(v)} $ 
where $ (X^{(u)},\tilde{X}^{(v)}) $ is considered 
under the product measure $ Q^{(u)} \otimes M^{(v)} $. 
Here and in what follows, we emphasize independence with the symbol $ \tilde{} $, 
unless otherwise stated. 

For $ t>0 $, we denote last exit time from the origin before $ t $ by 
\begin{align}
g_t = \inf \{ s \le t: X_s = 0 \} . 
\label{}
\end{align}
The following formula describes disintegration of $ P_0|_{\cF_t} $ 
with respect to $ g_t $: 

\begin{Thm} \label{thm: LED}
For each $ t>0 $, it holds that 
\begin{align}
P_0|_{\cF_t} = \int_0^t \vn(R>t-u) P_0[\d L_u] Q^{(u)} \bullet M^{(t-u)} . 
\label{eq: LED}
\end{align}
In other words, the following statements hold: 
\\ \quad 
{\rm (i)} 
The distribution of $ g_t $ is given by 
$ P_0(g_t \in \d u) = \vn(R>t-u) P_0[\d L_u] $; 
\\ \quad 
{\rm (ii)} 
Given $ g_t=u $, $ (X_t:t \in [0,u]) $ and $ (X_{u+t}:t \in [0,t-u]) $ 
are independent under $ P_0 $; 
\\ \quad 
{\rm (iii)} 
$ (X_t:t \in [0,u]) $ under $ P_0 $ is distributed as the stable bridge $ Q^{(u)} $; 
\\ \quad 
{\rm (iv)} 
$ (X_{u+t}:t \in [0,t-u]) $ under $ P_0 $ is distributed as the stable meander $ M^{(t-u)} $. 
\end{Thm}

\begin{Rem}
We note that the formula \eqref{eq: LED} is the counterpart 
of Salminen \cite[Prop. 4]{MR1454113} in his study of last exit decomposition 
for linear diffusions. 
\end{Rem}

\begin{Rem}
We remark that (i) implies 
\begin{align}
P_0(g_t \in \d u) = \frac{(t-u)^{\frac{1}{\alpha }-1} u^{-\frac{1}{\alpha }} \d u }
{\Gamma (1-1/\alpha ) \Gamma (1/\alpha )} 
\label{}
\end{align}
for some constant $ C $, which shows that 
$ \frac{1}{t} g_t $ has the Beta$ (1-\frac{1}{\alpha },\frac{1}{\alpha }) $ distribution. 
For further discussions, see \cite{YYY2}. 
\end{Rem}

\begin{proof}[Proof of Theorem \ref{thm: LED}]
Let us prove 
\begin{align}
P_0|_{\cF_t} 
= \int_0^t P_0[\d L_u] Q^{(u)} \bullet (\vn|_{\cF_{t-u}}) , 
\label{eq: P0(t)}
\end{align}
which is equivalent to \eqref{eq: LED}. 
Let $ F(t,w) $ be a non-negative continuous functional 
on $ [0,\infty ) \times \bD $. 
For each $ t \ge 0 $, we define a measurable functional $ F_t $ on $ \bD([0,t];\bR) $ 
by $ F_t(X^{(t)}) = F(t,X_{t \wedge \cdot}) $. 
Then we have 
\begin{align}
\int_0^{\infty } \d t F_t(X^{(t)}) 
= \sum_{l \in D} \int_0^{R(e_l)} \d r 
F_{\tau_{l-}+r} \cbra{ X^{(\tau_{l-})} \bullet e_l } . 
\label{}
\end{align}
Now we appeal to Theorem \ref{thm: compensation} 
and we obtain 
\begin{align}
\int_0^{\infty } P_0[F_t(X^{(t)})] \d t 
= 
\cbra{ P_0 \otimes \tilde{\vn} } 
\ebra{ \int_0^{\infty } \d L_t \int_0^{\infty } \d r 1_{ \{ \tilde{R}>r \} } 
F_{t+r} \cbra{ X^{(t)} \bullet \tilde{X}^{(r)} } } 
. 
\label{}
\end{align}
Since $ P_0[ \int_0^{\infty } G(X^{(u)}) \d L_u ] 
= \int_0^{\infty } P_0[\d L_u] Q^{(u)}[G(X^{(u)})] $, 
we obtain 
\begin{align}
\int_0^{\infty } P_0[F_t(X^{(t)})] \d t 
=& 
\int_0^{\infty } P_0[\d L_u] 
\cbra{ Q^{(u)} \otimes \tilde{\vn} } 
\ebra{ \int_0^{\infty } \d r 1_{ \{ \tilde{R}>r \} } 
F_{u+r} \cbra{ X^{(u)} \bullet \tilde{X}^{(r)} } } 
. 
\label{}
\end{align}
Changing variables to $ t=r+u $ and the order of integrations, we have 
\begin{align}
\int_0^{\infty } P_0[F_t(X^{(t)})] \d t 
= 
\int_0^{\infty } \d t 
\int_0^t P_0[\d L_u] 
\cbra{ Q^{(u)} \bullet (\vn|_{\cF_{t-u}}) } 
\ebra{ 1_{ \{ R>t-u \} } 
F_t \cbra{ X^{(t)} } } 
. 
\label{eq: int P0 Ft X(t)}
\end{align}
Since the identity \eqref{eq: int P0 Ft X(t)} holds 
with $ F_t $ replaced by $ \e^{-qt} F_t $ for any $ q>0 $, 
we obtain 
\begin{align}
P_0[F_t(X^{(t)})] 
= \int_0^t P_0[\d L_u] \cbra{ Q^{(u)} \bullet (\vn|_{\cF_{t-u}}) } 
\ebra{ 1_{ \{ R>t-u \} } F_t(X^{(t)}) } . 
\label{}
\end{align}
This completes the proof. 
\end{proof}

\begin{Rem}
In the above argument, we have proven the following formulae: 
\begin{align}
\int_0^{\infty } P_0^{(t)} \d t 
=& 
\int_0^{\infty } P_0^{(\tau_l)} \d l \bullet 
\int_0^{\infty } \vn(R>r) M^{(r)} \d r 
\\
=& 
\int_0^{\infty } P_0[\d L_u] Q^{(u)} \bullet 
\int_0^{\infty } \vn(R>r) M^{(r)} \d r . 
\label{}
\end{align}
Here we adopt the notations $ P_0^{(t)} $ and $ P_0^{(\tau_l)} $ 
which are found in \cite{MR1725357}, 
but we do not go into details. 
\end{Rem}

\subsection{Harmonicity of the function $ |x|^{\alpha -1} $}

Set 
\begin{align}
h(x) 
= \lim_{q \to 0+} \{ u_q(0)-u_q(x) \} 
= \frac{1}{\pi} \int_0^{\infty } \frac{1-\cos x \lambda}{\lambda^{\alpha }} \d \lambda . 
\label{}
\end{align}
Then we have 
\begin{align}
h(x) 
= h(1) |x|^{\alpha -1} 
\label{}
\end{align}
where 
\begin{align}
h(1) = 2 \cos \frac{(2-\alpha) \pi}{2} . 
\label{}
\end{align}

\begin{Thm} \label{thm: harmonic}
The function $ h(x) = h(1) |x|^{\alpha -1} $ is harmonic 
for the killed process, i.e., 
\begin{align}
P^0_x[h(X_t)] = P_x[h(X_t);T_{\{ 0 \}}>t] = h(x) 
, \qquad x \in \bR \setminus \{ 0 \} , \ t>0 . 
\label{}
\end{align}
Equivalently, $ (h(X_{t \wedge T_{\{ 0 \}}})) $ is a $ (P_x,\cF_t) $-martingale. 
\end{Thm}

We omit the proof, 
because Theorem \ref{thm: harmonic} follows immediately from the 

\begin{Thm}[Salminen--Yor \cite{SY}] \label{thm: Sal-Yor}
\Confirst{con: Tanaka}
For $ x \in \bR $, 
there exist a square-integrable martingale $ N_t^x $ 
and some constant $ C $ 
such that 
\begin{align}
|X_t|^{\alpha -1} = |x|^{\alpha -1} + N_t^x + C L(t,x) 
\quad \text{under} \ 
P_x . 
\label{}
\end{align}
\end{Thm}

\begin{Thm} \label{thm: harmonic2}
It holds that 
\begin{align}
\vn[h(X_t)] = 1 
, \qquad t>0 . 
\label{eq: vn hXt = 1}
\end{align}
\end{Thm}

\begin{proof}[Proof of Theorem \ref{thm: harmonic2}]
Theorem \ref{thm: rho} and the identity \eqref{eq: LT of T0} imply that 
\begin{align}
\int_0^{\infty } \e^{-qt} \vn[h(X_t)] \d t 
= \int h(x) P_x[\e^{-q T_{\{ 0 \}}}] \d x 
= \int h(x) \frac{u_q(x)}{u_q(0)} \d x . 
\label{}
\end{align}
Hence it suffices to prove that 
\begin{align}
\int u_q(x) h(x) \d x = \frac{u_q(0)}{q} 
, \qquad x \in \bR. 
\label{eq: uq h}
\end{align}
Let $ r $ be such that $ 0<r<q $. 
By the resolvent equation $ U_q U_r = (U_r-U_q)/(q-r) $, we have 
\begin{align}
\int u_q(x-y) u_r(y-z) \d y = \frac{1}{q-r} \kbra{ u_r(x-z) - u_q(x-z) } . 
\label{}
\end{align}
Letting $ x=z=0 $ and using the symmetry $ u_q(-y)=u_q(y) $, we have 
\begin{align}
\int u_q(y) u_r(y) \d y = \frac{1}{q-r} \kbra{ u_r(0) - u_q(0) } . 
\label{}
\end{align}
Now we have 
\begin{align}
\int u_q(y) \kbra{ u_r(0) - u_r(y) } \d y 
= \frac{u_q(0)}{q-r} - \frac{ru_r(0)}{q(q-r)} . 
\label{}
\end{align}
If we let $ r $ decrease to 0, then we see that 
\begin{align}
u_r(0) - u_r(x) 
= \frac{1}{\pi} \int_0^{\infty } \frac{1-\cos x \lambda}{r+\lambda^{\alpha }} \d \lambda 
\label{}
\end{align}
increases to $ h(x) $, 
and that $ ru_r(0) \to 0 $. Hence we obtain \eqref{eq: uq h} 
by the monotone convergence theorem. 
\end{proof}

\begin{Rem}
For generalisations of Theorems \ref{thm: harmonic} and \ref{thm: harmonic2} 
for symmetric L\'evy processes, see \cite{Y}. 
\end{Rem}

\subsection{Convergence of the stable meander to the $ h $-path process}

Let us introduce the {\em $ h $-path process} $ (P^{\dagger}_x:x \in \bR) $ 
as 
\begin{align}
P^{\dagger}_x |_{\cF_t} 
=& 
\frac{h(X_t)}{h(x)} \cdot P^0_x |_{\cF_t} 
, \qquad x \in \bR \setminus \{ 0 \}, 
\label{}
\\
P^{\dagger}_0 |_{\cF_t} 
=& 
h(X_t) \cdot \vn|_{\cF_t} . 
\label{eq: Imhof}
\end{align}
From Theorem \ref{thm: harmonic} and the Markov properties of $ P^0_x $ and $ \vn $, 
it follows that such a process exists uniquely. 
Remark that, when $ \alpha =2 $, 
the $ h $-path process coincides up to some scale transform 
with the symmetrization of three-dimensional Bessel process; 
consequently, the identity \eqref{eq: Imhof} 
is nothing but the {\em Imhof relation} 
(see, e.g., \cite[17, Exercise XII.4.18]{MR1725357}).

The following result asserts 
that the meander converges to the $ h $-path process. 

\begin{Thm} \label{thm: meander}
It holds that 
\begin{align}
M^{(t)} 
\ \stackrel{t \to \infty }{\longrightarrow} \ 
P^{\dagger}_0 
\qquad \text{along} \ 
(\cF_s) . 
\label{eq: meander conv Pdagger}
\end{align}
\end{Thm}

In order to prove Theorem \ref{thm: meander}, we need the 

\begin{Lem} \label{lem: Ytx}
For $ t>0 $ and $ x \neq 0 $, set 
\begin{align}
Y(t,x) = 
\frac{P_x(T_{\{ 0 \}}>t)}{h(x) \vn(R>t)} . 
\label{eq: Ytx}
\end{align}
Then it holds 
that $ Y(t,x) \to 1 $ as $ t \to \infty $ for any fixed $ x \neq 0 $, 
and that $ Y(t,x) $ is bounded in $ t>0 $ and $ x \neq 0 $. 
\end{Lem}

\begin{proof}[Proof of Lemma \ref{lem: Ytx}]
Using \eqref{eq: LT of T0}, we have 
\begin{align}
\int_0^{\infty } \e^{-qt} P_x(T_{\{ 0 \}}>t) \d t 
= \frac{u_q(0)-u_q(x)}{q u_q(0)} 
\sim h(x) \frac{q^{-1/\alpha } }{u_1(0)} 
\qquad \text{as} \ 
q \to 0+ . 
\label{}
\end{align}
Hence we may apply a tauberian theorem. 
By Proposition \ref{prop: vn Rt}, we obtain 
\begin{align}
P_x(T_{\{ 0 \}}>t) \sim 
h(x) \vn(R>t) 
\qquad \text{as $ t \to \infty $} . 
\label{asymp}
\end{align}
This shows the first assertion. 

Since the function $ t \mapsto Y(t,1) $ is continuous 
and $ Y(t,1) \to 1 $ as $ t \to \infty $, 
we see that $ Y(t,1) $ is bounded in $ t>0 $. 
By scaling property $ P_x(T_{\{ 0 \}}>t) = P_1(T_{\{ 0 \}}>|x|^{-\alpha }t) $, we have 
$ Y(t,x)= Y(|x|^{-\alpha }t,1) $. 
This proves the second assertion. 
\end{proof}

Now let us proceed to prove Theorem \ref{thm: meander}. 

\begin{proof}[Proof of Theorem \ref{thm: meander}]
Let $ s>0 $ be fixed 
and let $ Z_s $ be a bounded $ \cF_s $-measurable functional. 
By the Markov property of $ \vn $, we have 
\begin{align}
\vn \ebra{ Z_s 1_{\{ R>t \}} } 
= \vn \ebra{ Z_s 1_{\{ R>s \}} P_{X_s} (T_{\{ 0 \}}>t-s) } . 
\label{}
\end{align}
By the Imhof relation \eqref{eq: Imhof} and by \eqref{eq: Ytx}, we have 
\begin{align}
\vn \ebra{ Z_s 1_{\{ R>t \}} } 
=& P^{\dagger}_0 \ebra{ Z_s P_{X_s} (T_{\{ 0 \}}>t-s) / h(X_s) } 
\\
=& P^{\dagger}_0 [Z_s Y(t-s,X_s)] \cdot \vn(R>t-s) . 
\label{}
\end{align}
Dividing both sides by $ \vn(R>t) $, using Proposition \ref{prop: vn Rt}, 
and then applying the bounded convergence theorem, 
we obtain 
\begin{align}
M^{(t)}[Z_s] = P^{\dagger}_0[Z_s Y(t-s,X_s)] \cdot \cbra{ 1-\frac{s}{t} }^{\frac{1}{\alpha }-1} 
\to P^{\dagger}_0[Z_s] 
\label{eq: meander M and DP}
\end{align}
as $ t \to \infty $. 
This completes the proof. 
\end{proof}

\subsection{Convergence of the meander weighed by a multiplicative functional}

Let $ (\cE_t:t \ge 0) $ be an $ (\cF_t) $-adapted process 
which satisfies $ 0 \le \cE_t \le 1 $ 
and enjoys the multiplicativity property: 
\begin{align}
\cE_{t+s} = \cE_t \cdot (\cE_s \circ \theta_t) . 
\label{}
\end{align}
Such a process is called a {\em multiplicative functional}; 
see, e.g., \cite{MR0264757}. 
Then it necessarily follows that 
$ t \mapsto \cE_t $ is non-increasing. 

For later use, 
we need the following result 
which asserts that the convergence of the meander to the $ h $-path process 
is still valid with an extra weighing by a multiplicative functional. 

\begin{Thm} \label{thm: Penal P2 meander}
\begin{align}
\cE_t \cdot M^{(t)} 
\ \stackrel{t \to \infty }{\longrightarrow} \ 
\cE_{\infty } \cdot P^{\dagger}_0 
\qquad \text{along} \ 
(\cF_s) . 
\label{}
\end{align}
\end{Thm}

To prove Theorem \ref{thm: Penal P2 meander}, we need the following two lemmas.

\begin{Lem} \label{lem: conv in law}
For any $ x \in \bR $, it holds that 
\begin{align}
\begin{split}
& \cbra{ (X_t:t \ge 0),(\lambda^{-1/\alpha } X_{\lambda t}:t \ge 0) } 
\ \text{under} \ P^{\dagger}_x 
\\
\claw& 
\cbra{ (X_t:t \ge 0),(\tilde{X}_t:t \ge 0) } 
\ \text{under} \ P^{\dagger}_x \otimes \tilde{P} 
\end{split}
\label{eq: conv of pair in law}
\end{align}
as $ \lambda \to \infty $ 
where 
\begin{align}
\tilde{P} = 
\begin{cases}
P^{\dagger}_0 \ & \text{if} \ 1<\alpha <2 \ \text{or if} \ x=0, 
\\
P^{\rm 3B,+}_0 \ & \text{if} \ \alpha =2 \ \text{and} \ x>0, 
\\
P^{\rm 3B,-}_0 \ & \text{if} \ \alpha =2 \ \text{and} \ x<0. 
\end{cases}
\label{}
\end{align}
\end{Lem}

\begin{proof}
We prove the claim only in the case $ 1<\alpha <2 $; 
in fact, almost the same argument works in the other cases. 
Set $ X^{(\lambda)}_t = \lambda^{-1/\alpha } X_{\lambda t} $. 
Let us apply the convergence theorem of \cite[Theorem VI.16]{MR762984}. 

First, let $ t \ge 0 $ be fixed 
and let $ f:\bR \to \bR $ be a continuous function such that 
$ \lim_{|x| \to \infty } f(x) = 0 $. 
Then we have 
\begin{align}
\lim_{\lambda \to \infty } P^{\dagger}_x[f(X^{(\lambda)}_t)] 
= 
\lim_{\lambda \to \infty } P^{\dagger}_{\lambda^{-1/\alpha } x}[f(X_t)] 
= 
P^{\dagger}_0[f(X_t)] . 
\label{}
\end{align}
In fact, the first identity follows from the scaling property 
and the second follows from the Feller property of the $ h $-path process, 
which is proved in \cite{Y}. 
Hence we obtain 
\begin{align}
X^{(\lambda)}_t \ \text{under} \ P^{\dagger}_x 
\quad \claw \quad 
X_t \ \text{under} \ P^{\dagger}_0 
\label{}
\end{align}
as $ \lambda \to \infty $. 
By a standard argument involving the Markov property, we see that 
the convergence \eqref{eq: conv of pair in law} holds 
in the sense of finite dimensional distributions. 

Second, for any sequence $ \{ \lambda_n \} $ with $ \lambda_n \to \infty $, 
let us check the {\em Aldous condition}: 
For a sequence of positive constants $ \{ \delta_n \} $ converging to zero 
and for a bounded sequence of stopping times $ \{ \rho_n \} $, 
\begin{align}
\abra{ X_{\rho_n+\delta_n} - X_{\rho_n} } 
+ 
\abra{ X^{(\lambda_n)}_{\rho_n+\delta_n} - X^{(\lambda_n)}_{\rho_n} } 
\stackrel{n \to \infty }{\longrightarrow} 
0 
\qquad \text{in $ P^{\dagger}_x $-probability}. 
\label{eq: Aldous}
\end{align}
The convergence \eqref{eq: Aldous} is equivalent to 
\begin{align}
X^{(\lambda_n)}_{\rho_n+\delta_n} 
- X^{(\lambda_n)}_{\rho_n} 
\stackrel{n \to \infty }{\longrightarrow} 
0 
\qquad \text{in $ P^{\dagger}_x $-probability}. 
\label{eq: Aldous2}
\end{align}
To prove \eqref{eq: Aldous2}, it suffices to prove that 
\begin{align}
P^{\dagger}_x \ebra{ \abra{ X^{(\lambda_n)}_{\rho_n+\delta_n} 
- X^{(\lambda_n)}_{\rho_n} } \wedge 1 } 
\stackrel{n \to \infty }{\longrightarrow} 
0 . 
\label{eq: Aldous3}
\end{align}
By the strong Markov property and by the scaling property, we have 
\begin{align}
P^{\dagger}_x \ebra{ \abra{ X^{(\lambda_n)}_{\rho_n+\delta_n} 
- X^{(\lambda_n)}_{\rho_n} } \wedge 1 } 
= 
P^{\dagger}_x \ebra{ 
P^{\dagger}_{\lambda^{-1/\alpha } X_{\rho_n}} 
\ebra{ \abra{ X_{\delta_n} - X_0 } \wedge 1 } } . 
\label{}
\end{align}
Hence we can easily obtain the convergence \eqref{eq: Aldous3} 
by the Feller property of the $ h $-path process. 
\end{proof}

\begin{Lem} \label{lem: Penal P2 meander}
For any $ x \neq 0 $, it holds that 
\begin{align}
\frac{P_x[\cE_t;T_{\{ 0 \}}>t]}{h(x) \vn(R>t)} 
\to P^{\dagger}_x[\cE_{\infty }] 
\qquad \text{as} \ 
t \to \infty . 
\label{}
\end{align}
\end{Lem}

\begin{proof}[Proof of Lemma \ref{lem: Penal P2 meander}]
For $ t>s>0 $, we have $ \cE_t \le \cE_s $, 
and hence we have 
\begin{align}
\frac{P_x[\cE_t;T_{\{ 0 \}}>t]}{h(x) \vn(R>t)} 
\le& 
\frac{P_x[\cE_s;T_{\{ 0 \}}>t]}{h(x) \vn(R>t)} 
\\
=& 
P^{\dagger}_x \ebra{ \cE_s \frac{P_{X_s}(T_{\{ 0 \}}>t-s)}{h(X_s) \vn(R>t)} } 
\\
=& 
P^{\dagger}_x [ \cE_s Y(t-s,X_s) ] \cdot \cbra{1-\frac{s}{t}}^{\frac{1}{\alpha }-1} . 
\label{}
\end{align}
By Lemma \ref{lem: Ytx} and by the bounded convergence theorem, 
we have 
\begin{align}
\limsup_{t \to \infty } 
\frac{P_x[\cE_t;T_{\{ 0 \}}>t]}{h(x) \vn(R>t)} 
\le P^{\dagger}_x[\cE_s] . 
\label{}
\end{align}
Since $ P^{\dagger}_x[\cE_s] \to P^{\dagger}_x[\cE_{\infty }] $ as $ s \to \infty $, 
we obtain the upper estimate: 
\begin{align}
\limsup_{t \to \infty } 
\frac{P_x[\cE_t;T_{\{ 0 \}}>t]}{h(x) \vn(R>t)} 
\le P^{\dagger}_x[\cE_{\infty }] . 
\label{}
\end{align}

By Lemma \ref{lem: conv in law}, we have 
\begin{align}
\cbra{ \cE_t,t^{-1/\alpha } |X_t| } 
\ \text{under} \ P^{\dagger}_x 
\claw 
(\cE_{\infty },|\tilde{X}_1|) 
\ \text{under} \ P^{\dagger}_x \otimes P^{\dagger}_0 . 
\label{}
\end{align}
Hence, by Fatou's lemma, we have 
\begin{align}
\liminf_{t \to \infty } 
\frac{P^{\dagger}_x[\cE_t / h(X_t)]}{\vn(R>t)} 
=& 
\liminf_{t \to \infty } 
\frac{P^{\dagger}_x[\cE_t / h(t^{-1/\alpha } |X_t|)]}{\vn(R>1)} 
\\
\ge& 
\frac{P^{\dagger}_x[\cE_{\infty }] P^{\dagger}_0[1/h(|X_1|)]}{\vn(R>1)} 
= 
P^{\dagger}_x[\cE_{\infty }] . 
\label{}
\end{align}
Thus we obtain the lower estimate: 
\begin{align}
\liminf_{t \to \infty } 
\frac{P_x[\cE_t;T_{\{ 0 \}}>t]}{h(x) \vn(R>t)} 
= \liminf_{t \to \infty } 
\frac{P^{\dagger}_x[\cE_t / h(X_t)]}{\vn(R>t)} 
\ge P^{\dagger}_x[\cE_{\infty }] . 
\label{}
\end{align}
Therefore the proof is now completed. 
\end{proof}

Now we prove Theorem \ref{thm: Penal P2 meander}. 

\begin{proof}[Proof of Theorem \ref{thm: Penal P2 meander}]
For a bounded $ \cF_s $-measurable functional $ Z_s $ and for $ t>s>0 $, we have 
\begin{align}
M^{(t)}[Z_s \cE_t] 
= P^{\dagger}_0 \ebra{ Z_s \cE_s \frac{P_{X_s}[\cE_{t-s};T_{\{ 0 \}}>t-s]}{h(X_s) \vn(R>t-s)} } 
\cdot \cbra{ 1-\frac{s}{t} }^{\frac{1}{\alpha }-1} . 
\label{}
\end{align}
Note that 
\begin{align}
\frac{P_x[\cE_t;T_{\{ 0 \}}>r]}{h(x) \vn(R>r)} 
\le 
\frac{P_x(T_{\{ 0 \}}>r)}{h(x) \vn(R>r)} 
= Y(r,x) , 
\label{}
\end{align}
which is uniformly bounded in $ r>0 $ and $ x \neq 0 $ by Lemma \ref{lem: Ytx}. 
Note also that 
\begin{align}
\frac{P_x[\cE_t;T_{\{ 0 \}}>r]}{h(x) \vn(R>r)} 
\ \stackrel{r \to \infty }{\longrightarrow} \ 
P^{\dagger}_x[\cE_{\infty }] 
, \qquad x \neq 0 
\label{}
\end{align}
by Lemma \ref{lem: Penal P2 meander}. 
Hence we apply bounded convergence theorem and obtain 
\begin{align}
M^{(t)}[Z_s \cE_t] 
\ \stackrel{t \to \infty }{\longrightarrow} \ 
P^{\dagger}_0 \ebra{ Z_s \cE_s P^{\dagger}_{X_s}[\cE_{\infty }] } 
= 
P^{\dagger}_0 \ebra{ Z_s \cE_{\infty } } . 
\label{}
\end{align}
This completes the proof. 
\end{proof}

\section{General observations on the $ \sigma $-finite measure 
unifying our penalisation problems 
and the martingale generator}
\label{sec: univ}

Following \cite{NRY} and \cite{NRY2}, we make general observations 
on the measure $ \sP $. 

\subsection{The $ \sigma $-finite measure unifying our penalisation problems}

Recall the definition of $ \sP $: 
\begin{align}
\sP 
= \int_0^{\infty } P_0[\d L_u] Q^{(u)} \bullet P^{\dagger}_0 
\label{eq: cP}
\end{align}
where 
\begin{align}
P_0[\d L_u] = \frac{\Gamma(1/\alpha )}{\alpha \pi} \frac{\d u}{u^{1/\alpha }} 
\label{}
\end{align}
and where 
$ P^{\dagger}_0 $ is defined by 
\begin{align}
P^{\dagger}_0 |_{\cF_t} = h(X_t) \cdot \vn |_{\cF_t} 
, \qquad t>0 . 
\label{}
\end{align}
Denote 
\begin{align}
g = \sup \{ t \ge 0 : X_t=0 \} . 
\label{}
\end{align}

\begin{Thm} \label{thm: sP}
The following statements hold: 
\\ \quad 
{\rm (i)} 
$ \sP(g \in \d u) = P_0[\d L_u] $; 
\\ \quad 
{\rm (ii)} 
$ \sP $ is a $ \sigma $-finite measure on $ \cF_{\infty } $; 
\\ \quad 
{\rm (iii)} 
$ \sP $ is singular with respect to $ P_0 $ on $ \cF_{\infty } $; 
\\ \quad 
{\rm (iv)} 
For each $ t>0 $ and for $ A \in \cF_t $, 
one has 
\begin{align}
\sP(A) =& 0 \qquad \text{if} \ P_0(A)=0, 
\label{eq: sPA = 0}
\\
\sP(A) =& \infty \qquad \text{if} \ P_0(A)>0. 
\label{eq: sPA = infty}
\end{align}
\end{Thm}

\begin{Rem}
For each $ t>0 $, \eqref{eq: sPA = 0} asserts that 
$ \sP $ is equivalent to $ P_0 $ on $ \cF_t $, 
but \eqref{eq: sPA = infty} asserts that 
$ \sP $ is {\em never} $ \sigma $-finite on $ \cF_t $. 
We insist that, 
since $ \sP $ is not $ \sigma $-finite on $ \cF_t $, 
\eqref{eq: sPA = 0} does not imply the existence 
of an $ \cF_t $-measurable Radon--Nikodym density. 
\end{Rem}

\begin{proof}[Proof of Theorem \ref{thm: sP}]
{\rm (i)} 
Since $ P^{\dagger}_0 $ is locally equivalent to $ \vn $, 
we see that 
$ P^{\dagger}_0(X_s \neq 0 \ \text{for any} \ s \le t) = 1 $ 
for any $ t>0 $. 
This shows that 
$ P^{\dagger}_0(X_t \neq 0 \ \text{for any} \ t>0) = 1 $. 
Hence we see, by the definition \eqref{eq: cP} of $ \sP $, that 
$ g=u $ under the measure $ Q^{(u)} \bullet P^{\dagger}_0 $. 
Thus we obtain the desired result. 

{\rm (ii)} 
It is obvious by (i) that 
$ \sP(g<u) $ is finite for each $ u>0 $. 

{\rm (iii)} 
On one hand, we have $ \sP(g=\infty )=0 $. 
On the other hand, 
since the origin for $ (X_t,P_0) $ is recurrent, we have $ P_0(g<\infty )=0 $. 
This implies that $ \sP $ is singular to $ P_0 $ on $ \cF_{\infty } $. 

{\rm (iv)} 
Let $ A \in \cF_t $ 
and suppose that $ P_0(A)=0 $. 
For $ T>t $, we have 
\begin{align}
\int_0^T P_0[\d L_u] \cbra{ Q^{(u)} \bullet P^{\dagger}_0 }(A) 
= 
\int_0^t P_0[\d L_u] \cbra{ Q^{(u)} \bullet P^{\dagger}_0 }(A) 
+ \int_t^T P_0[\d L_u] Q^{(u)}(A) . 
\label{}
\end{align}
For $ 0<u<t $, we have 
$ (Q^{(u)} \bullet P^{\dagger}_0)(A) 
= (Q^{(u)} \bullet \vn)[1_A h(X_t)] $, 
and hence we obtain 
\begin{align}
\int_0^t P_0[\d L_u] \cbra{ Q^{(u)} \bullet P^{\dagger}_0 }(A) = P_0 \ebra{ 1_A h(X_t) } = 0 . 
\label{}
\end{align}
For $ t<u<T $, we have 
\begin{align}
\int_t^T P_0[\d L_u] Q^{(u)}(A) = P_0 \ebra{ 1_A (L_T-L_t) } = 0 . 
\label{}
\end{align}
Letting $ T \to \infty $, we obtain $ \int_t^{\infty } P_0[\d L_u] Q^{(u)}(A) = 0 $. 
Therefore we obtain $ \sP(A)=0 $. 

Conversely, 
let $ A \in \cF_t $ 
and suppose that $ P_0(A)>0 $. 
Then 
\begin{align}
\sP(A) \ge \int_t^{\infty } P_0[\d L_u] Q^{(u)}(A) = P_0 \ebra{ 1_A (L_{\infty }-L_t) } . 
\label{}
\end{align}
Note that the last quantity is $ \infty $ 
since $ P_0(L_{\infty }=\infty )=1 $. 
Hence we obtain $ \sP(A)=\infty $. 
\end{proof}

\subsection{The martingale generator}

\begin{Thm} \label{thm: mart op}
For each $ x \in \bR $, $ t \ge 0 $ 
and for each non-negative measurable or $ \sP_x $-integrable functional $ F $, 
there exists a unique $ \cF_t $-measurable functional $ M_{t,x}(F) $ (possibly taking infinite values) such that 
\begin{align}
(F \cdot \sP_x)|_{\cF_t} = M_{t,x}(F) \cdot P_x|_{\cF_t} . 
\label{}
\end{align}
In particular, if $ F $ is $ \sP_x $-integrable, 
then the process $ (M_{t,x}(F):t \ge 0) $ is an $ (\cF_t,P_x) $-martingale 
such that 
\begin{align}
M_{0,x}(F)=\sP_x[F] 
\label{}
\end{align}
and that 
\begin{align}
\lim_{t \to \infty } M_{t,x}(F) = 0 
\qquad \text{$ P_x $-almost surely.}
\label{}
\end{align}
\end{Thm}

In the case $ x=0 $, we write $ M_t(F) $ for $ M_{t,0}(F) $. 
For each $ x \in \bR $, 
we call the operator $ L^1(\sP_x) \ni F \mapsto (M_{t,x}(F):t \ge 0) $ 
the {\em martingale generator}. 

\begin{proof}
It is obvious that the uniqueness holds in the sense that, 
if $ F=G $ $ \sP $-almost everywhere, then $ M_{t,x}(F)=M_{t,x}(G) $ $ P_0 $-almost surely. 
Without loss of generality, 
we may suppose that $ x=0 $ and that $ F $ is non-negative. 

Let $ n $ be a positive integer and set $ F_n = F \cdot 1_{\{ g<n \}} $. 
By (ii) and (iv) of Theorem \ref{thm: sP}, 
we see that $ (F_n \cdot \sP)|_{\cF_t} $ is 
a finite measure and is absolutely continuous with respect to $ P_0|_{\cF_t} $. 
Hence we may apply the Radon--Nikodym theorem 
to obtain the desired functional $ M_t(F_n) $ as the Radon--Nikodym derivative. 
Hence the desired functional $ M_t(F) $ is obtained 
as the increasing limit $ \lim_{n \to \infty } M_t(F_n) $ 
by the monotone convergence theorem. 

Suppose that $ F $ is $ \sP $-integrable. 
For $ s \le t $, we have 
\begin{align}
P_0[Z_s M_t(F)] 
= \sP[Z_s F] 
= P_0[Z_s M_s(F)] . 
\label{}
\end{align}
Hence $ (M_t(F):t \ge 0) $ is a $ (\cF_t,P_0) $-martingale. 
It is obvious that 
$ M_0(F)=\sP[F] $. 

Since $ (M_t(F):t \ge 0) $ is a non-negative martingale, 
$ M_t(F) $ converges $ P_0 $-almost surely 
to a non-negative $ \cF_{\infty } $-measurable functional $ M_{\infty }(F) $. 
For $ 0<s<t \le \infty $, set $ A(s,t) = \{ g_t \ge s \} \in \cF_t $. 
Note that $ P_0(A(s,\infty )) = P_0(g \ge s) = 1 $. 
Applying Fatou's lemma 
and then applying the dominated convergence theorem, 
we obtain 
\begin{align}
P_0[M_{\infty }(F)] 
= P_0[1_{A(s,\infty )} M_{\infty }(F)] 
\le& 
\liminf_{t \to \infty } P_0[1_{A(s,t)} M_t(F)] 
\\
=& 
\liminf_{t \to \infty } \sP[1_{A(s,t)} F] 
= \sP[1_{A(s,\infty )} F] . 
\label{}
\end{align}
Since $ \sP(g=\infty )=0 $, we have $ \lim_{s \to \infty } \sP[1_{A(s,\infty )} F] = 0 $. 
Hence we obtain $ P_0[M_{\infty }(F)]=0 $, 
which implies that $ P_0(M_{\infty }(F)=0)=1 $. 
Therefore the proof is completed. 
\end{proof}

\section{Convergence lemmas}
\label{sec: imp lem}

Let $ 0<\gamma <1 $. 
For integrable functions $ \psi_t(u) $ 
such that $ \psi_t(u) \to \exists \psi(u) $ as $ t \to \infty $, 
we may expect that 
\begin{align}
\int_0^t \cbra{1-\frac{u}{t}}^{\gamma -1} \psi_t(u) \d u 
\to 
\int_0^{\infty } \psi(u) \d u 
\qquad \text{as} \ t \to \infty . 
\label{eq: imp lem asymp}
\end{align}
We need this convergence for several functions $ \psi_t $ 
in order to solve our penalisation problems, 
as we have seen roughly in \eqref{eq: convergence of P_0/nRt to sP}. 
In fact, we shall see that 
we must be careful in dealing with the convergence \eqref{eq: imp lem asymp}. 
In this section 
we give some sufficient conditions for the convergence \eqref{eq: imp lem asymp} 
as well as a counterexample. 

If $ \psi_t $'s satisfy 
\begin{align}
\int_0^t \psi_t(u) \d u \to \int_0^{\infty } \psi(u) \d u 
\qquad \text{as} \ t \to \infty , 
\label{eq: imp lem asymp2}
\end{align}
then the convergence \eqref{eq: imp lem asymp} is equivalent to 
\begin{align}
I(\psi_t,t) \to 0 
\qquad \text{as} \ 
t \to \infty 
\label{}
\end{align}
where 
\begin{align}
I(\psi,t) = \int_0^t \kbra{ \cbra{ 1-\frac{u}{t} }^{\gamma -1} -1 } \psi(u) \d u . 
\label{}
\end{align}

First, we present the following counterexample. 

\begin{Ex}
The convergence \eqref{eq: imp lem asymp} {\em fails} if 
\begin{align}
\psi_t(u) \equiv \psi(u) 
= \sum_{n=1}^{\infty } n^{\frac{2+\gamma }{1-\gamma }} 
1_{\cbra{n-n^{-\frac{4-\gamma }{1-\gamma }},n}}(u) . 
\label{}
\end{align}
\end{Ex}

\begin{proof}
$ \psi $ is integrable since 
$ \int_0^{\infty } \psi(u) \d u = \sum_{n=1}^{\infty } n^{-2} < \infty $. 
But $ \limsup_t I(\psi,t) = \infty $ because 
\begin{align}
I(\psi,n) 
\ge& n^{\frac{2+\gamma }{1-\gamma }} \cdot n^{1-\gamma } 
\int_{n-n^{-\frac{4-\gamma }{1-\gamma }}}^n (n-u)^{\gamma -1} \d u 
- n^{-2} 
\\
=& n^{\frac{2+\gamma }{1-\gamma }} \cdot n^{1-\gamma } 
\cdot \gamma ^{-1} n^{-\frac{(4-\gamma )\gamma }{1-\gamma }} - n^{-2} 
\\
=& \gamma ^{-1} n^{3-2\gamma } - n^{-2} 
\to \infty 
\qquad \text{as} \ 
n \to \infty . 
\label{}
\end{align}
This prevents the convergence \eqref{eq: imp lem asymp}. 
\end{proof}

On the other hand, we give three sufficient conditions 
for the convergence \eqref{eq: imp lem asymp}; 
the first one is rather theoretical, 
but the second and third ones can be readily applied. 

\begin{Lem}[Dominated convergence] \label{lem: imp lem1}
Suppose that $ \psi_t $'s are integrable functions 
such that $ \int_0^{\infty } \psi_t(u) \d u \to \int_0^{\infty } \psi(u) \d u $ 
for some integrable function $ \psi $. 
Suppose, in addition, that $ |\psi_t| \le \tilde{\psi}_t $ 
for some integrable function $ \tilde{\psi}_t $ 
such that $ \lim_{t \to \infty } I(\tilde{\psi}_t,t) = 0 $. 
Then 
\begin{align}
\int_0^t \cbra{1-\frac{u}{t}}^{\gamma -1} \psi_t(u) \d u 
\to 
\int_0^{\infty } \psi(u) \d u 
\qquad \text{as} \ t \to \infty 
\label{}
\end{align}
holds. 
\end{Lem}

\begin{proof}
This is obvious 
by $ |I(\psi_t,t)| \le I(|\psi_t|,t) \le I(\tilde{\psi}_t,t) \to 0 $ as $ t \to \infty $. 
\end{proof}

\begin{Lem} \label{lem: imp lem3}
Suppose that $ \psi $ is a non-negative integrable function and satisfies 
\begin{align}
\lim_{t \to \infty } \kbra{ t \sup_{u>t} \psi(u) } = 0 . 
\label{eq: assump sup psi}
\end{align}
Then $ \lim_{t \to \infty } I(\psi,t) = 0 $. 
\end{Lem}

\begin{proof}
Let $ 0<\eps<1 $ be fixed. 
We split $ I(\psi,t) $ into a sum $ I(\psi_1,t) + I(\psi_2,t) $ 
where 
$ \psi_1 = \psi 1_{(\eps t,\infty )} $ 
and 
$ \psi_2 = \psi 1_{(0,\eps t)} $. 

By the definition of $ I(\psi_1,t) $ and changing variables to $ v=ut $, we have 
\begin{align}
I(\psi_1,t) 
=& 
\int_{\eps t}^t 
\kbra{ \cbra{ 1-\frac{u}{t} }^{\gamma -1} -1 } \psi(u) \d u 
\\
\le& 
t \sup_{u>\eps t} \psi(u) 
\int_{\eps t}^t 
\kbra{ \cbra{ 1-\frac{u}{t} }^{\gamma -1} -1 } \frac{\d u}{t} 
\\
=& 
\frac{1}{\eps} \kbra{ \eps t \sup_{u>\eps t} \psi(u) } 
\int_{\eps}^1 \kbra{ (1-v)^{\gamma -1} - 1 } \d v . 
\label{}
\end{align}
By the assumption \eqref{eq: assump sup psi}, we obtain 
$ \lim_{t \to \infty } I(\psi_1,t) = 0 $ 
for any fixed $ \eps>0 $. 

By the definition of $ I(\psi_2,t) $, we have 
\begin{align}
I(\psi_2,t) 
\le 
\kbra{ \frac{1}{(1-\eps)^{1-\gamma }} -1 } 
\int_0^{\infty } \psi(u) \d u . 
\label{}
\end{align}
Hence we have $ \limsup_{t \to \infty } I(\psi_2,t) $ vanishes as $ \eps \to 0+ $. 
Now the proof is completed. 
\end{proof}

\begin{Lem} \label{lem: imp lem2}
Suppose that $ \psi_t(u)=\psi_1(u) \psi_2(t-u) $ 
where $ \psi_1 $ is integrable 
and $ \psi_2 $ is bounded measurable 
with $ \lim_{u \to \infty } \psi_2(u) = \psi_2(\infty ) > 0 $. 
Suppose, in addition, that the function 
$ t \mapsto \int_0^t (t-u)^{\gamma -1} \psi_t(u) \d u $ 
is ultimately non-increasing as $ t $ increases. 
Then $ \lim_{t \to \infty } I(\psi_t,t) = 0 $. 
\end{Lem}

\begin{proof}
Taking the Laplace transform, we have 
\begin{align}
\int_0^{\infty } \d t \e^{-qt} \int_0^t (t-u)^{\gamma -1} \psi_t(u) \d u 
=& 
\int_0^{\infty } \e^{-qu} \psi_1(u) \d u 
\int_0^{\infty } \e^{-qt} t^{\gamma -1} \psi_2(t) \d t 
\\
\sim& 
\Gamma (\gamma ) q^{-\gamma } \psi_2(\infty ) \int_0^{\infty } \psi_1(u) \d u 
\qquad \text{as} \ 
q \to 0+ . 
\label{}
\end{align}
Hence we may apply the tauberian theorem. 
By the monotonicity assumption, we obtain 
\begin{align}
\int_0^t (t-u)^{\gamma -1} \psi_t(u) \d u 
\sim 
t^{\gamma -1} \psi_2(\infty ) \int_0^{\infty } \psi_1(u) \d u 
\qquad \text{as} \ t \to \infty . 
\label{}
\end{align}
On the other hand, we have 
\begin{align}
\int_0^t \psi_t(u) \d u 
= \int_0^t \psi_1(u) \psi_2(t-u) \d u 
\to 
\psi_2(\infty ) \int_0^{\infty } \psi_1(u) \d u 
\qquad \text{as} \ t \to \infty . 
\label{}
\end{align}
Therefore we obtain $ \lim_{t \to \infty } I(\psi_t,t) = 0 $. 
\end{proof}

\section{Penalisation with a function of the local time at the origin}
\label{sec: LT penal}

\subsection{Results}

\begin{Thm} \label{thm: LT mart op}
Let $ f $ be a non-negative function on $ [0,\infty ) $. Then 
it holds that 
\begin{align}
M_t(f(L_{\infty })) 
= h(X_t) f(L_t) + \int_{L_t}^{\infty } f(l) \d l 
, \qquad t \ge 0 . 
\label{eq: LT mart op}
\end{align}
Consequently, it holds that 
\begin{align}
\sP[f(L_{\infty })] 
= M_0(f(L_{\infty })) 
= \int_0^{\infty } f(l) \d l . 
\label{}
\end{align}
\end{Thm}

\begin{Rem}
As an outcome of \eqref{eq: LT mart op}, 
we have established that its right hand side is a $ (P_0,\cF_t) $-martingale, 
a well-known fact for $ \alpha =2 $ (see \cite[Prop. VI.4.5]{MR1725357})
\end{Rem}

\begin{Thm} \label{thm: LT Penal}
Let $ f $ be a non-negative function on $ [0,\infty ) $ such that 
\begin{align}
I(f) := \int_0^{\infty } f(l) \d l \in (0,\infty ) . 
\label{eq: LT Penal assump}
\end{align}
Then it holds that 
\begin{align}
\frac{f(L_t) \cdot P_0}{\vn(R>t)} 
\ \stackrel{t \to \infty }{\longrightarrow} \ 
f(L_{\infty }) \cdot \sP 
\qquad \text{along $ (\cF_s) $.} 
\label{eq: LT conv}
\end{align}
Consequently, the penalisation with the weight functional $ \Gamma_t=f(L_t) $ 
is given as 
\begin{align}
\frac{f(L_t) \cdot P_0}{P_0[f(L_t)]} 
\ \stackrel{t \to \infty }{\longrightarrow} \ 
\frac{f(L_{\infty }) \cdot \sP}{I(f)} 
\qquad \text{along $ (\cF_s) $.} 
\label{}
\end{align}
\end{Thm}

\subsection{Proofs}

\begin{proof}[Proof of Theorem \ref{thm: LT mart op}]
Let $ t>0 $ be fixed 
and $ Z_t $ a non-negative $ \cF_t $-measurable functional. 
On the one hand, 
since $ L_{\infty } = L_t $ on $ \{ g \le t \} $, we have 
\begin{align}
\sP[Z_t f(L_{\infty }) 1_{\{ g \le t \}}] 
=& \sP[Z_t f(L_t) 1_{\{ g \le t \}}] 
\\
=& \int_0^t P_0[\d L_u] \cbra{ Q^{(u)} \bullet P^{\dagger}_0 } [Z_t f(L_t)] 
\\
=& \int_0^t P_0[\d L_u] \cbra{ Q^{(u)} \bullet \vn } [Z_t f(L_t) h(X_t)] 
\\
=& \int_0^t \vn(R>t-u) P_0[\d L_u] \cbra{ Q^{(u)} \bullet M^{(t-u)} } [Z_t f(L_t) h(X_t)] 
\\
=& P_0 [Z_t f(L_t) h(X_t)] . 
\label{}
\end{align}
On the other hand, we have 
\begin{align}
\sP[Z_t f(L_{\infty }) 1_{\{ g>t \}}] 
=& \int_t^{\infty } P_0[\d L_u] \cbra{ Q^{(u)} \bullet P^{\dagger}_0 } [Z_t f(L_u)] 
\\
=& \int_t^{\infty } P_0[\d L_u] Q^{(u)} [Z_t f(L_u)] 
\\
=& P_0 \ebra{ Z_t \int_t^{\infty } f(L_u) \d L_u } 
\\
=& P_0 \ebra{ Z_t \int_{L_t}^{\infty } f(l) \d l } . 
\label{}
\end{align}
Hence we obtain 
\begin{align}
\sP[Z_t f(L_{\infty })] 
= P_0 \ebra{ Z_t \kbra{ f(L_t) h(X_t) + \int_{L_t}^{\infty } f(l) \d l } } . 
\label{}
\end{align}
Therefore we have completed the proof. 
\end{proof}

\begin{proof}[Proof of Theorem \ref{thm: LT Penal}]
We need only to prove the first assertion that 
\begin{align}
\frac{f(L_t) \cdot P_0}{\vn(R>t)} 
\ \stackrel{t \to \infty }{\longrightarrow} \ 
f(L_{\infty }) \cdot \sP 
\qquad \text{along $ (\cF_s) $.} 
\label{}
\end{align}

Set $ \psi(u) = p_u(0) Q^{(u)}[f(L_u)] $. 
We will prove in Lemma \ref{lem: Penal P1 psi} below 
that $ \psi $ satisfies the assumption of Lemma \ref{lem: imp lem3}. 
Now we apply Lemma \ref{lem: imp lem3} for the function $ \psi $ 
and we obtain 
\begin{align}
\int_0^t & \cbra{1-\frac{u}{t}}^{\frac{1}{\alpha }-1} P_0[\d L_u] Q^{(u)}[f(L_u)] 
\to 
\int_0^{\infty } P_0[\d L_u] Q^{(u)}[f(L_u)] 
\label{eq: penal conv P1-0}
\end{align}
as $ t \to \infty $. 
Let $ s>0 $ be fixed 
and let $ Z_s $ be a bounded $ \cF_s $-measurable functional. 
Then 
\begin{align}
\int_0^t P_0[\d L_u] \cbra{ Q^{(u)} \bullet M^{(t-u)} } [Z_s f(L_u)] 
\to 
\int_0^{\infty } P_0[\d L_u] \cbra{ Q^{(u)} \bullet P^{\dagger}_0 } [Z_s f(L_u)] 
\label{}
\end{align}
as $ t \to \infty $ 
by Lebesgue's convergence theorem. 
Hence we can apply Lemma \ref{lem: imp lem1}, and we obtain 
\begin{align}
\frac{P_0[Z_s f(L_t)]}{\vn(R>t)} 
=& \int_0^t \cbra{1-\frac{u}{t}}^{\frac{1}{\alpha }-1} 
P_0[\d L_u] \cbra{ Q^{(u)} \bullet M^{(t-u)} } [Z_s f(L_u)] 
\\
\to&
\int_0^{\infty } P_0[\d L_u] \cbra{ Q^{(u)} \bullet P^{\dagger}_0 } [Z_s f(L_u)] 
\qquad \text{(as $ t \to \infty $)} 
\label{eq: penal conv P1}
\\
=& \sP[Z_s f(L_{\infty })] 
\end{align}
as $ t \to \infty $. 
This completes the proof. 
\end{proof}

\begin{Lem} \label{lem: Penal P1 psi}
Set 
\begin{align}
\psi(u) = p_u(0) Q^{(u)}[f(L_u)] . 
\label{eq: Penal P1 psi}
\end{align}
Then the function $ \psi(u) $ is continuous and $ u \psi(u) \to 0 $ as $ u \to \infty $. 
In particular, the function $ \psi $ satisfies 
the assumption \eqref{eq: assump sup psi} of Lemma \ref{lem: imp lem3}. 
\end{Lem}

\begin{proof}
For any non-negative Borel function $ \phi $, we have 
\begin{align}
\int_0^{\infty } \phi(u) \psi(u) \d u 
=& P_0 \ebra{ \int_0^{\infty } \phi(u) f(L_u) \d L_u } 
\\
=& \int_0^{\infty } P_0[\phi(\tau_s)] f(s) \d s 
\\
=& \int_0^{\infty } P_0[\phi(s^{1/\beta } \tau_1)] f(s) \d s 
\intertext{
where $ \beta = 1 - 1/\alpha $. 
If we denote $ \rho^{(\beta )}(v) = P_0(\tau_1 \in \d v) / \d v $, we have }
=& \int_0^{\infty } \d s f(s) \int_0^{\infty } \phi(s^{1/\beta } v) \rho^{(\beta )}(v) \d v 
\\
=& \int_0^{\infty } \d u \phi(u) \int_0^{\infty } s^{-1/\beta } 
\rho^{(\beta )}(s^{-1/\beta } u) f(s) \d s . 
\label{}
\end{align}
Hence we obtain 
\begin{align}
\psi(u) = 
\int_0^{\infty } s^{-1/\beta } \rho^{(\beta )}(s^{-1/\beta } u) f(s) \d s . 
\label{}
\end{align}
Since the function $ \rho^{(\beta )}(v) $ is unimodal (see, e.g., Sato \cite{MR1739520}), 
we see that 
$ v \rho^{(\beta )}(v) $ is bounded in $ v>0 $ 
and that $ v \rho^{(\beta )}(v) \to 0 $ as $ v \to \infty $. 
Therefore, by the assumption that $ \int_0^{\infty } f(s) \d s < \infty $, 
we obtain the desired result. 
\end{proof}

\begin{Rem}
In the Brownian case $ \alpha =2 $, 
the corresponding $ \beta $ equals $ 1/2 $ and 
\begin{align}
\rho^{(1/2)}(v) = \frac{1}{2 \sqrt{\pi v^3}} \e^{-\frac{1}{4v}} . 
\label{}
\end{align}
\end{Rem}

\section{Feynman--Kac penalisations}
\label{sec: FK penal}

\subsection{Results}

Recall that our Feynman--Kac penalisation is the penalisation with the weight functional 
\begin{align}
\cE^V_t = \exp \kbra{ - \int L(t,x) V(\d x) } 
, \qquad t \ge 0 
\label{}
\end{align}
for a non-negative measure $ V(\d x) $ on $ \bR $.

\begin{Thm} \label{thm: FK Penal}
Let $ V $ be a non-negative measure on $ \bR $ such that 
\begin{align}
0 < \int (1+|y|^{\alpha -1}) V(\d y) < \infty . 
\label{eq: FK Penal assump}
\end{align}
Let $ x \in \bR $. Then it holds that 
\begin{align}
0 < \sP_x[\cE^V_{\infty }] < \infty 
\label{}
\end{align}
and that 
\begin{align}
\frac{( \cE^V_t 1_{\{ T_{\{ 0 \}} > t \}} ) \cdot P_x}{\vn(R>t)} 
\ \stackrel{t \to \infty }{\longrightarrow}& \ 
( \cE^V_{\infty } 1_{\{ T_{\{ 0 \}} = \infty \}} ) \cdot \sP_x 
\qquad \text{along $ (\cF_s) $,} 
\label{eq: FK conv 1}
\\
\frac{( \cE^V_t 1_{\{ T_{\{ 0 \}} \le t \}} ) \cdot P_x}{\vn(R>t)} 
\ \stackrel{t \to \infty }{\longrightarrow}& \ 
( \cE^V_{\infty } 1_{\{ T_{\{ 0 \}} < \infty \}} ) \cdot \sP_x 
\qquad \text{along $ (\cF_s) $} 
\label{eq: FK conv 2}
\intertext{and}
\frac{\cE^V_t \cdot P_x}{\vn(R>t)} 
\ \stackrel{t \to \infty }{\longrightarrow}& \ 
\cE^V_{\infty } \cdot \sP_x 
\qquad \text{along $ (\cF_s) $.} 
\label{eq: FK conv}
\end{align}
\end{Thm}

\begin{Cor}
Let $ V $ be a non-negative measure on $ \bR $ such that \eqref{eq: FK Penal assump} holds. 
Then the penalisation with the weight functional 
$ \Gamma_t = \cE^V_t $ is given as 
\begin{align}
\frac{\cE^V_t \cdot P_x}{P_x[\cE^V_t]} 
\ \stackrel{t \to \infty }{\longrightarrow} \ 
\frac{\cE^V_{\infty } \cdot \sP_x}{\sP_x[\cE^V_{\infty }]} 
\qquad \text{along $ (\cF_s) $.} 
\label{eq: FK penal}
\end{align}
\end{Cor}

\begin{Thm} \label{thm: FK mart op}
Let $ V $ be a non-negative measure on $ \bR $ such that \eqref{eq: FK Penal assump} holds. 
Set 
\begin{align}
C_V = \sP[\cE^V_{\infty }] . 
\label{}
\end{align}
Let $ x \in \bR $. 
Then it holds that 
\begin{align}
\varphi^1_V(x) 
:=& \lim_{t \to \infty } \frac{P_x[\cE^V_t;T_{\{ 0 \}}>t]}{\vn(R>t)} 
= \sP_x[\cE^V_{\infty };T_{\{ 0 \}}=\infty ] 
= h(x) P^{\dagger}_x[\cE^V_{\infty }] , 
\label{eq: FK func 1}
\\
\varphi^2_V(x) 
:=& \lim_{t \to \infty } \frac{P_x[\cE^V_t;T_{\{ 0 \}} \le t]}{\vn(R>t)} 
= \sP_x[\cE^V_{\infty };T_{\{ 0 \}}<\infty ] 
= C_V P_x[\cE^V_{T_{\{ 0 \}}}] 
\label{eq: FK func 2}
\intertext{and} 
\varphi_V(x) 
:=& \lim_{t \to \infty } \frac{P_x[\cE^V_t]}{\vn(R>t)} 
= \sP_x[\cE^V_{\infty }] 
= h(x) P^{\dagger}_x[\cE^V_{\infty }] + C_V P_x[\cE^V_{T_{\{ 0 \}}}] 
\\
\equiv & 
\varphi^1_V(x) + \varphi^2_V(x) . 
\label{eq: FK func}
\end{align}
Moreover, for $ t \ge 0 $, it holds that 
\begin{align}
M_{t,x}(\cE^V_{\infty } 1_{\{ T_{\{ 0 \}} = \infty \}}) 
=& \varphi^1_V(X_t) 1_{\{ T_{\{ 0 \}} > t \}} \cE^V_t , 
\label{eq: FK mart op 1}
\\
M_{t,x}(\cE^V_{\infty } 1_{\{ T_{\{ 0 \}} < \infty \}}) 
=& \kbra{ 
\varphi^1_V(X_t) 1_{\{ T_{\{ 0 \}} \le t \}} 
+ \varphi^2_V(X_t) 
} \cE^V_t , 
\label{eq: FK mart op 2}
\intertext{and} 
M_{t,x}(\cE^V_{\infty }) 
=& \varphi_V(X_t) \cE^V_t . 
\label{eq: FK mart op}
\end{align}
\end{Thm}

We divide the proofs of Theorems \ref{thm: FK Penal} and \ref{thm: FK mart op} 
into several steps in the following subsections.

\begin{Rem}\footnote{
Tildes in this remark have nothing to do with our previous notation's 
for independence.} 
For the Feynman--Kac penalisations (Theorems \ref{thm: FK Penal} and \ref{thm: FK mart op}) 
in the brownian case, 
Roynette--Vallois--Yor (\cite{MR2261065},\cite{MR2229621}, \cite{MR2253307}) 
have given more characterisations of the limit measure 
than the contents of Theorem \ref{thm: FK mart op}. 
For convenience, we consider the Wiener measures $ (W_x:x \in \bR) $ 
normalized with the weight functional 
\begin{align}
\tilde{\cE}^V_t = \exp \kbra{ - \frac{1}{2} \int L(t,x) V(\d x) } . 
\label{}
\end{align}
For each $ x \in \bR $, let $ \sW_x $ denote the law of $ (x+X_t:t \ge 0) $ under $ \sW $. 
Then the function 
\begin{align}
\tilde{\varphi}_V(x) 
= \lim_{t \to \infty } \sqrt{t} W_x[\tilde{\cE}^V_t] 
= \sW_x[\tilde{\cE}^V_{\infty }] 
\label{}
\end{align}
is the unique solution of the Sturm--Liouville differential equation 
\begin{align}
\d \tilde{\varphi}'_V(x) = \tilde{\varphi}_V(x) V(\d x) 
\label{eq: Sturm--Liouville}
\end{align}
subject to the boundary conditions 
\begin{align}
\lim_{x \to \infty } \tilde{\varphi}_V'(x) = \sqrt{\frac{2}{\pi}} 
\qquad \text{and} \qquad 
\lim_{x \to - \infty } \tilde{\varphi}_V'(x) = - \sqrt{\frac{2}{\pi}} . 
\label{}
\end{align}
Moreover, the limit measure $ W^V_x $ (instead of $ P^V_x $) 
is the law of the unique solution of the stochastic differential equation 
\begin{align}
\d X_t = \d B_t + \frac{\tilde{\varphi}'_V(X_t)}{\tilde{\varphi}_V(X_t)} \d t 
, \qquad X_0 = x . 
\label{eq: SDE}
\end{align}
We do not know how to develop these arguments in the stable L\'evy case, 
for which it would be interesting to obtain counterparts of 
\eqref{eq: Sturm--Liouville} and \eqref{eq: SDE}. 
\end{Rem}

\subsection{Penalisation weighed by a general multiplicative functional}

In this subsection, we make a general study. 
Let $ \cE = (\cE_t:t \ge 0) $ be an $ (\cF_t) $-adapted process 
which satisfies $ 0 \le \cE_t \le 1 $, $ t \ge 0 $ 
and 
is a multiplicative functional: 
\begin{align}
\cE_{t+s} = \cE_t \cdot (\cE_s \circ \theta_t) 
, \qquad t,s \ge 0 . 
\label{}
\end{align}
Note that the process $ t \mapsto \cE_t $ is necessarily non-increasing; 
in fact, $ \cE_{t+s} = \cE_t \cdot (\cE_s \circ \theta_t) \le \cE_t $ 
for any $ t,s \ge 0 $. 

\begin{Thm} \label{thm: general FK penal}
Let $ x \in \bR $ be fixed. 
Suppose that 
\begin{align}
\int_0^{\infty } P_0[\d L_u] Q^{(u)} [\cE_u] < \infty . 
\label{eq: assump for general FK penal}
\end{align}
Then it holds that 
\begin{align}
\frac{\cE_t \cdot P_x}{\vn(R>t)} 
\ \stackrel{t \to \infty }{\longrightarrow} \ 
\cE_{\infty } \cdot \sP_x 
\qquad \text{along} \ 
(\cF_s) . 
\label{}
\end{align}
\end{Thm}

\begin{proof}
We may suppose that $ x=0 $ without loss of generality. 

By the multiplicativity property, we have 
\begin{align}
\sP[\cE_{\infty }] 
=& \int_0^{\infty } P_0[\d L_u] \cbra{ Q^{(u)} \bullet P^{\dagger}_0 } [\cE_{\infty }] 
\\
=& \kbra{ \int_0^{\infty } P_0[\d L_u] Q^{(u)}[\cE_u] } P^{\dagger}_0[\cE_{\infty }] . 
\label{}
\end{align}

Set 
\begin{align}
\psi_t(u) = p_u(0) (Q^{(u)} \bullet M^{(t-u)}) [\cE_t] . 
\label{}
\end{align}
Then we have $ \psi_t(u) = \psi_1(u) \psi_2(t-u) $ 
where 
$ \psi_1(u) = p_u(0) Q^{(u)}[\cE_u] $ 
and $ \psi_2(t) = M^{(t)}[\cE_t] $. 
Let us check that 
all the assumptions of Lemma \ref{lem: imp lem2} are satisfied for $ \psi_t(u) $. 
Note that 
\begin{align}
\int_0^{\infty } \psi_1(u) \d u 
= \int_0^{\infty } P_0[\d L_u] Q^{(u)} [\cE_u] 
\label{}
\end{align}
and it is finite by the assumption \eqref{eq: assump for general FK penal}. 
Note also that $ \psi_2 $ is bounded and that 
$ \lim_{t \to \infty } \psi_2(t) = P^{\dagger}_0[\cE_{\infty }] $ 
by Theorem \ref{thm: Penal P2 meander}. 
Recall the following identity: 
\begin{align}
P_0[\cE_t] 
= \int_0^t \vn(R>t-u) P_0[\d L_u] \cbra{ Q^{(u)} \bullet M^{(t-u)} } [\cE_t] . 
\label{}
\end{align}
Since the left-hand side is non-increasing as $ t $ increases, 
we see that the function $ t \mapsto \int_0^t (t-u)^{\frac{1}{\alpha }-1} \psi_t(u) \d u $ 
is non-increasing as $ t $ increases. 
Hence we have verified all the assumptions of Lemma \ref{lem: imp lem2}, 
and we obtain $ \lim_{t \to \infty } I(\psi_t,t) = 0 $. 
The remainder of the proof follows from Lemma \ref{lem: imp lem1}. 
\end{proof}

\subsection{Non-degeneracy condition}

Now we return to the case where $ \cE_t = \cE^V_t $. 
By the multiplicativity of $ (\cE^V_t) $, we have 
\begin{align}
C_V 
= \sP[\cE^V_{\infty }] 
=& \int_0^{\infty } P_0[\d L_u] \cbra{ Q^{(u)} \bullet P^{\dagger}_0 } [\cE^V_{\infty }] 
\\
=& \int_0^{\infty } P_0[\d L_u] Q^{(u)}[\cE^V_u] P^{\dagger}_0[\cE^V_{\infty }] 
\\
=& \kbra{ \int_0^{\infty } P_0[\cE^V_{\tau_s}] \d s } P^{\dagger}_0[\cE^V_{\infty }] . 
\label{}
\end{align}

\begin{Thm} \label{thm: conti case}
The following assertions hold: 
\\ \quad 
{\rm (i)} 
If $ V \neq 0 $, then $ \displaystyle \int_0^{\infty } P_0[\cE^V_{\tau_s}] \d s < \infty $; 
\\ \quad 
{\rm (ii)} 
If $ V((-\eps,\eps)) < \infty $ for some $ \eps>0 $, 
then $ \displaystyle \int_0^{\infty } P_0[\cE^V_{\tau_s}] \d s > 0 $; 
\\ \quad 
{\rm (iii)} 
If $ \displaystyle \int h(x) V(\d x) < \infty $, 
then $ P^{\dagger}_0[\cE^V_{\infty }] > 0 $; 
\\ \quad 
{\rm (iv)} 
If $ \displaystyle 0 < \int \{ 1+h(x) \} V(\d x) < \infty $, 
then $ 0<C_V<\infty $. 
\end{Thm}

For the proof of Theorem \ref{thm: conti case}, we need the following 

\begin{Lem} \label{lem: expect of LT}
The following statements hold: 
\\ \quad {\rm (i)} 
$ \vn[L(R,x)] = 1 $ for any $ x \in \bR \setminus \{ 0 \} $; 
\\ \quad {\rm (ii)} 
$ P^{\dagger}_0[L(t,x)] = h(x) P_x(T_{\{ 0 \}}<t) $ for any $ t \ge 0 $ and any $ x \in \bR $; 
\\ \quad {\rm (iii)} 
$ P^{\dagger}_0[L(\infty ,x)] = h(x) $ for any $ x \in \bR $. 
\end{Lem}

Remark that $ \vn[L(R,0)]=0 $; 
in fact, $ L(R,0)=0 $ $ \vn $-almost everywhere. 

\begin{proof}
(i) 
For a non-negative Borel function $ f $, we have 
\begin{align}
\int \d x f(x) \vn[L(R,x)] 
= \int_0^{\infty } \vn[f(X_t)] \d t 
= \int \d x f(x) \int_0^{\infty } \rho(t,x) \d t 
= \int \d x f(x) . 
\label{}
\end{align}
Hence we obtain $ \vn[L(R,x)] = 1 $ for almost every $ x \in \bR $. 
By the scaling property, we obtain the desired conclusion. 

(ii) 
Let $ t \ge 0 $ be fixed. 
For a non-negative Borel function $ f $, we have 
\begin{align}
\int \d x f(x) P^{\dagger}_0[L(t,x)] 
=& \int_0^t P^{\dagger}_0[f(X_s)] \d s 
= \int_0^t \vn[f(X_s)h(X_s)] \d s 
\\
=& \int \d x f(x) h(x) \int_0^t \rho(s,x) \d s 
\\
=& \int \d x f(x) h(x) P_x(T_{\{ 0 \}}<t) . 
\label{}
\end{align}
Hence we see that 
$ P^{\dagger}_0[L(t,x)] = h(x) P_x(T_{\{ 0 \}}<t) $ for almost every $ x \in \bR $. 
Since $ t \mapsto P^{\dagger}_0[L(t,1)] $ is continuous by the monotone convergence theorem, 
we see, by the scaling property, 
that $ \bR \setminus \{ 0 \} \ni x \mapsto P^{\dagger}_0[L(t,x)] $ 
is continuous. 
Noting that $ L(t,0)=0 $ $ P^{\dagger}_0 $-almost surely, 
we complete the proof. 

(iii) 
Letting $ t \to \infty $ in (ii), 
we obtain 
$ P^{\dagger}_0[L(\infty ,x)] = h(x) $ 
by the monotone convergence theorem. 
\end{proof}

Now we prove Theorem \ref{thm: conti case}. 

\begin{proof}[Proof of Theorem \ref{thm: conti case}]
Note that 
\begin{align}
\cE^V_{\tau_s} 
= \exp \kbra{ - s V(\{ 0 \}) 
- \sum_{l \in D, \ l \le s} \int_{\{ x \neq 0 \}} L(R,x)[e_l] V(\d x) } 
\label{}
\end{align}
where $ L(R,x)[e_l] $ is the local time at $ x $ of the excursion $ e_l $ up to its lifetime. 
Hence we have 
$ P_0[ \cE^V_{\tau_s} ] = \exp \kbra{ - s K_V } $ 
where 
\begin{align}
K_V = V(\{ 0 \}) + \vn \ebra{ 1-\exp \kbra{ - \int_{\{ x \neq 0 \}} L(R,x) V(\d x) } } . 
\label{}
\end{align}
Consequently we have 
\begin{align}
\int_0^{\infty } P_0 \ebra{ \cE_{\tau_s} } \d s = \frac{1}{K_V} . 
\label{}
\end{align}

{\rm (i)} 
If $ \int_0^{\infty } P_0[\cE^V_{\tau_s}] \d s = \infty $, 
then we have $ K_V=0 $, which implies that $ V=0 $. 
Hence the assertion is proved by contraposition. 

{\rm (ii)} 
Suppose that $ V((-\eps,\eps))<\infty $ for $ \eps>0 $. 
Then, by Lemma \ref{lem: expect of LT}, we have 
\begin{align}
\vn \ebra{ \int_{(-\eps,\eps)} L(R,x) V(\d x) } = \int_{(-\eps,\eps)} \vn[ L(R,x) ] V(\d x) 
= V((-\eps,\eps))<\infty . 
\label{}
\end{align}
Now we obtain 
\begin{align}
& \vn \ebra{ 1-\exp \kbra{ - \int L(R,x) V(\d x) } ; \sup_{t \ge 0} |X(t)| <\eps } 
\\
=& \vn \ebra{ 1-\exp \kbra{ - \int_{(-\eps,\eps)} L(R,x) V(\d x) } ; \sup_{t \ge 0} |X(t)| <\eps } 
\\
\le& \vn \ebra{ \int_{(-\eps,\eps)} L(R,x) V(\d x) } < \infty . 
\label{}
\end{align}
Since $ \vn(\sup_{t \ge 0} |X(t)| \ge \eps) < \infty $, we obtain $ K_V<\infty $. 
Hence the assertion is proved. 

{\rm (iii)} 
By Lemma \ref{lem: expect of LT}, we obtain 
\begin{align}
P^{\dagger}_0 \ebra{ \int L(\infty ,x) V(\d x) } 
= \int h(x) V(\d x) 
< \infty . 
\label{}
\end{align}
This implies that 
\begin{align}
P^{\dagger}_0 \cbra{ \int L(\infty ,x) V(\d x) < \infty } 
= P^{\dagger}_0(\cE^V_{\infty }>0) 
= 1 , 
\label{}
\end{align}
which proves $ P^{\dagger}_0[\cE^V_{\infty }]>0 $. 

{\rm (iv)} 
Suppose that $ 0 < \int \{ 1+h(x) \} V(\d x) < \infty $. 
Then the assumptions of (i)-(iii) are all satisfied. 
Noting that $ \cE^V_{\infty } \le 1 $, 
we obtain $ 0 < C_V < \infty $. 
\end{proof}

\subsection{Proof of Theorems}

\begin{proof}[Proof of Theorem \ref{thm: FK Penal}]
Note that $ (\cE^V_t) $ and $ (\cE^V_t 1_{\{ T_{\{ 0 \}} > t \}}) $ 
are multiplicative functionals which take values in $ [0,1] $. 
By Theorem \ref{thm: conti case}, 
we may apply Theorem \ref{thm: general FK penal} 
to obtain \eqref{eq: FK conv} and \eqref{eq: FK conv 1}. 
Subtracting both sides of \eqref{eq: FK conv 1} from \eqref{eq: FK conv}, 
we obtain \eqref{eq: FK conv 2}. 
\end{proof}

\begin{proof}[Proof of Theorem \ref{thm: FK mart op}]
The second equalities of \eqref{eq: FK func 1}, \eqref{eq: FK func 2} and \eqref{eq: FK func} 
are obvious by Theorem \ref{thm: FK Penal}. 
The last equality of \eqref{eq: FK func 1} is obvious 
by Lemma \ref{lem: Penal P2 meander}. 
The last equality of \eqref{eq: FK func 2} is obtained as follows: 
\begin{align}
\frac{P_x[\cE^V_t 1_{\{ T_{\{ 0 \}} \le t \}}]}{\vn(R>t)} 
=& 
\int_0^t P_x[\cE^V_s;T_{\{ 0 \}} \in \d s] \frac{P_0[\cE^V_{t-s}]}{\vn(R>t)} 
\\
\stackrel{t \to \infty }{\longrightarrow}& 
\int_0^{\infty } P_x[\cE^V_s;T_{\{ 0 \}} \in \d s] \sP[\cE^V_{\infty }] 
= C_V P_x[\cE^V_{T_{\{ 0 \}}}] . 
\label{}
\end{align}
Now we obtain the last equality of \eqref{eq: FK func} 
by adding \eqref{eq: FK func 1} and \eqref{eq: FK func 2}. 

Let $ 0 \le s < t $. 
By the Markov property of $ (P_x) $ and by Theorem \ref{thm: FK Penal}, we have 
\begin{align}
P_x \ebra{ \cE^V_s 1_{\{ T_{\{ 0 \}} >s \}} 
\frac{ P_{X_s}[\cE^V_{t-s};T_{\{ 0 \}} >t-s ]}{\vn(R>t)} } 
=& 
\frac{P_x[\cE^V_t;T_{\{ 0 \}} >t ]}{\vn(R>t)} 
\\
\stackrel{t \to \infty }{\longrightarrow}& 
\sP_x[\cE^V_{\infty } 1_{\{ T_{\{ 0 \}} = \infty \}} ] 
\\
=& 
h(x) P^{\dagger}_x[\cE^V_{\infty }] . 
\label{}
\end{align}
By the Markov property of $ (P^{\dagger}_x) $, we have 
\begin{align}
h(x) P^{\dagger}_x[\cE^V_{\infty }] 
= h(x) P^{\dagger}_x[\cE^V_s P^{\dagger}_{X_s}[\cE^V_{\infty }] ] 
= P_x[\cE^V_s 1_{\{ T_{\{ 0 \}} >s \}} \varphi^1_V(X_s)] . 
\label{}
\end{align}
Hence, by Scheff\'e's lemma, we obtain 
\begin{align}
\frac{ P_{X_s}[\cE^V_{t-s};T_{\{ 0 \}} >t-s ]}{\vn(R>t)} 
\stackrel{t \to \infty }{\longrightarrow}& 
\varphi^1_V(X_s) 
\qquad \text{in $ L^1(\cE^V_s 1_{\{ T_{\{ 0 \}} >s \}} \cdot P_x) $.} 
\label{}
\end{align}
Therefore, for any bounded $ \cF_s $-measurable functional $ Z_s $, we have 
\begin{align}
\frac{P_x[Z_s \cE^V_t;T_{\{ 0 \}} >t ]}{\vn(R>t)} 
=& 
P_x \ebra{ Z_s \cE^V_s 1_{\{ T_{\{ 0 \}} >s \}} 
\frac{ P_{X_s}[\cE^V_{t-s};T_{\{ 0 \}} >t-s ]}{\vn(R>t)} } 
\label{eq: FK penal conv1}
\\
\stackrel{t \to \infty }{\longrightarrow}& 
P_x [ Z_s \cE^V_s 1_{\{ T_{\{ 0 \}} >s \}} \varphi^1_V(X_s) ] . 
\label{}
\end{align}
Combining this with \eqref{eq: FK conv 1}, we obtain 
\begin{align}
\sP_x[Z_s \cE^V_{\infty } 1_{\{ T_{\{ 0 \}} = \infty \}} ] 
= 
P_x [ Z_s \cE^V_s 1_{\{ T_{\{ 0 \}} >s \}} \varphi^1_V(X_s) ] . 
\label{}
\end{align}
This implies the identity \eqref{eq: FK mart op 1}. 

By similar arguments, we have 
\begin{align}
\frac{P_x[Z_s \cE^V_t 1_{\{ T_{\{ 0 \}} \le s \}} ]}{\vn(R>t)} 
=& P_x \ebra{ Z_s \cE^V_s 1_{\{ T_{\{ 0 \}} \le s \}} \frac{P_{X_s}[\cE^V_{t-s}]}{\vn(R>t)} } 
\\
\stackrel{t \to \infty }{\longrightarrow}& 
P_x \ebra{ Z_s \cE^V_s 1_{\{ T_{\{ 0 \}} \le s \}} \varphi_V(X_s) } 
\label{}
\end{align}
and 
\begin{align}
\frac{P_x[Z_s \cE^V_t 1_{\{ s < T_{\{ 0 \}} \le t \}} ]}{\vn(R>t)} 
=& P_x \ebra{ Z_s \cE^V_s 1_{\{ T_{\{ 0 \}} > s \}} 
\frac{P_{X_s}[\cE^V_{t-s};T_{\{ 0 \}} \le t-s]}{\vn(R>t)} } 
\\
\stackrel{t \to \infty }{\longrightarrow}& 
P_x \ebra{ Z_s \cE^V_s 1_{\{ T_{\{ 0 \}} > s \}} \varphi^2_V(X_s) } . 
\label{}
\end{align}
Combining these two limits together with \eqref{eq: FK conv 2}, 
we obtain \eqref{eq: FK mart op 2}. 

The remainder of the proof is now obvious. 
\end{proof}

\section{Characterisation of non-negative martingales} \label{sec: further}

For a non-negative $ \sP $-integrable functional $ G $ such that $ \sP[G]>0 $, 
we define the probability measure $ P^G $ on $ \cF_{\infty } $ as 
\begin{align}
P^G = \frac{G \cdot \sP}{\sP[G]} . 
\label{}
\end{align}
We say that a statement holds {\em $ \sP $-almost surely} 
if it holds $ P^G $-almost surely 
for some $ \sP $-integrable functional $ G $ such that 
$ G>0 $ $ \sP $-almost everywhere. 
By the Radon--Nikodym theorem, 
$ \sP $-almost sure statement does not depend on the particular choice of such a functional $ G $. 

The following theorem is the stable L\'evy version of \cite[Corollary 1.2.6]{NRY}. 

\begin{Thm} \label{thm: char of pos mart}
Let $ (N_t) $ be a non-negative $ (\cF_t,P_0) $-martingale. 
Then $ (N_t) $ is represented as $ N_t = M_t(F) $ 
for some $ F \in L^1(\sP) $ 
if and only if it holds that 
\begin{align}
\frac{N_t}{1+h(X_t)} 
\stackrel{t \to \infty }{\longrightarrow} 
F 
\ \text{$ \sP $-almost surely} 
\quad \text{and} \quad 
\sP[F] = N_0 . 
\label{eq: conv sP a.s.}
\end{align}
\end{Thm}

Although it is completely parallel to that of \cite{NRY}, 
we give the proof for completeness of the paper.

\begin{Lem} \label{lem: ratio mart conv}
Let $ F $ and $ G $ be a non-negative $ \sP $-integrable functional 
and suppose that $ G>0 $ $ \sP $-almost everywhere. 
Then it holds that 
\begin{align}
\frac{M_t(F)}{M_t(G)} = P^G \ebra{ \frac{F}{G} \biggm| \cF_t } . 
\label{}
\end{align}
Consequently, it holds that 
\begin{align}
\frac{M_t(F)}{M_t(G)} 
\stackrel{t \to \infty }{\longrightarrow} 
\frac{F}{G} 
\qquad \text{$ P^G $-almost surely.} 
\label{}
\end{align}
\end{Lem}

\begin{proof}
Let $ Z_t $ be a non-negative $ \cF_t $-measurable functional. 
On the one hand, we have 
\begin{align}
P^G[Z_t F/G] = \sP[Z_t F] = P_0[Z_t M_t(F)] . 
\label{}
\end{align}
On the other hand, we have 
\begin{align}
P^G[Z_t M_t(F)/M_t(G)] = \sP[Z_t (M_t(F)/M_t(G)) G] 
= P_0[Z_t M_t(F)] . 
\label{}
\end{align}
Hence we obtain $ P^G[Z_t F/G]=P^G[Z_t M_t(F)/M_t(G)] $, which completes the proof. 
\end{proof}

\begin{Lem} \label{lem: MtF/1+hXt}
Let $ F $ be a non-negative $ \sP $-integrable functional. 
Then 
\begin{align}
\frac{M_t(F)}{1+h(X_t)} 
\stackrel{t \to \infty }{\longrightarrow} 
F 
\qquad \text{$ \sP $-almost surely.} 
\label{}
\end{align}
\end{Lem}

\begin{proof}
We apply Theorem \ref{thm: LT mart op} with $ f(l)=\e^{-l} $ to see that 
$ G=\e^{-L_{\infty }} $ is a positive $ \sP $-integrable functional such that 
\begin{align}
M_t(G) = (1+h(X_t)) \e^{-L_t} . 
\label{}
\end{align}
Hence we obtain 
\begin{align}
\frac{M_t(G)}{1+h(X_t)} 
\stackrel{t \to \infty }{\longrightarrow} 
G 
\qquad \text{$ P^G $-almost surely.} 
\label{}
\end{align}
Hence, by Lemma \ref{lem: ratio mart conv}, we obtain 
\begin{align}
\frac{M_t(F)}{1+h(X_t)} 
= 
\frac{M_t(G)}{1+h(X_t)} 
\cdot 
\frac{M_t(F)}{M_t(G)} 
\stackrel{t \to \infty }{\longrightarrow} 
G \cdot \frac{F}{G} 
= F 
\qquad \text{$ P^G $-almost surely.} 
\label{}
\end{align}
This completes the proof. 
\end{proof}

The following proposition states an interesting representation 
of any non-negative $ (P_0,\cF_t) $-supermartingale, 
a component of which is a certain $ (M_t(F)) $ martingale.

\begin{Prop} \label{prop: supermart decomp}
Let $ (N_t) $ a non-negative $ (\cF_t,P_0) $-supermartingale. 
\\ \quad {\rm (i)} 
There exists a non-negative $ \sP $-integrable functional $ F $ such that 
\begin{align}
\frac{N_t}{1+h(X_t)} \stackrel{t \to \infty }{\longrightarrow} F 
\qquad \text{$ \sP $-almost surely;} 
\label{}
\end{align}
\quad {\rm (ii)} 
Denote the $ P_0 $-almost sure limit of $ (N_t) $ as $ t \to \infty $ by $ N_{\infty } $. 
Then $ (N_t) $ decomposes uniquely in the following form: 
\begin{align}
N_t = M_t(F) + P_0[N_{\infty }|\cF_t] + \xi_t 
\label{eq: supermart decomp}
\end{align}
where: 
\\ \quad {\rm (iia)} 
$ (M_t(F)) $ is a non-negative $ (\cF_t,P_0) $-martingale such that 
\begin{align}
M_t(F) \stackrel{t \to \infty }{\longrightarrow} 0 
\ \text{($ P_0 $-a.s.)} 
\qquad \text{and} \qquad 
\frac{M_t(F)}{1+h(X_t)} \stackrel{t \to \infty }{\longrightarrow} F 
\ \text{($ \sP $-a.s.);} 
\label{}
\end{align}
\quad {\rm (iib)} 
$ (P_0[N_{\infty }|\cF_t]) $ is a non-negative uniformly-integrable $ (\cF_t,P_0) $-martingale 
with $ P_0 $-integrable terminal value $ N_{\infty } $ such that 
\begin{align}
P_0[N_{\infty }|\cF_t] \stackrel{t \to \infty }{\longrightarrow} N_{\infty } 
\ \text{($ P_0 $-a.s.)} 
\qquad \text{and} \qquad 
\frac{P_0[N_{\infty }|\cF_t]}{1+h(X_t)} \stackrel{t \to \infty }{\longrightarrow} 0 
\ \text{($ \sP $-a.s.);} 
\label{}
\end{align}
\quad {\rm (iic)} 
$ (\xi_t) $ is a non-negative $ (\cF_t,P_0) $-supermartingale such that 
\begin{align}
\xi_t \stackrel{t \to \infty }{\longrightarrow} 0 
\ \text{($ P_0 $-a.s.)} 
\qquad \text{and} \qquad 
\frac{\xi_t}{1+h(X_t)} \stackrel{t \to \infty }{\longrightarrow} 0 
\ \text{($ \sP $-a.s.)} 
\label{}
\end{align}
\end{Prop}

\begin{proof}
(i) 
Let $ G=\e^{-L_{\infty }} $. 
For any non-negative $ \cF_s $-measurable functional $ Z_s $, 
we see that 
\begin{align}
P^G[Z_s N_t/M_t(G)] = \sP[Z_s N_t G/M_t(G)] = P_0[Z_s N_t] \le P_0[Z_s N_s] . 
\label{eq: Nt/MtG is supermart}
\end{align}
Hence we conclude that 
$ (N_t/M_t(G)) $ is a non-negative $ (\cF_t,P^G) $-supermartingale. 
Thus there exists a non-negative $ \cF_{\infty } $-measurable functional $ \zeta $ such that 
$ N_t/M_t(G) \to \zeta $ $ P^G $-almost surely. 
By Lemma \ref{lem: MtF/1+hXt}, we see that 
\begin{align}
\frac{N_t}{1+h(X_t)} 
= \frac{N_t}{M_t(G)} 
\cdot \frac{M_t(G)}{1+h(X_t)} 
\stackrel{t \to \infty }{\longrightarrow} 
\zeta G =: F 
\qquad \text{$ P^G $-almost surely.} 
\end{align}

(ii) 
For any non-negative $ \cF_t $-measurable functional $ Z_t $, we have 
\begin{align}
P_0[Z_t M_t(F)] = \sP[Z_t F] = P^G[Z_t \zeta] . 
\end{align}
By Fatou's lemma, the last expectation is dominated by 
\begin{align}
\liminf_{u \to \infty } P^G \ebra{ Z_t \cdot \frac{N_u}{M_u(G)} } 
= \liminf_{u \to \infty } P_0[Z_t N_u] \le P_0[Z_t N_t] . 
\end{align}
This proves that $ M_t(F) \le N_t $ $ P_0 $-almost surely. 
Now we see that $ (\bar{N}_t := N_t - M_t(F)) $ 
is a non-negative $ (\cF_t,P_0) $-supermartingale. 
Since $ M_t(G) \to 0 $ $ P_0 $-almost surely as $ t \to \infty $, 
we see that 
\begin{align}
\lim_{t \to \infty } \bar{N}_t 
= 
\lim_{t \to \infty } N_t 
= 
N_{\infty } 
\qquad \text{$ P_0 $-almost surely.} 
\end{align}
For any non-negative $ \cF_t $-measurable functional $ Z_t $, we have 
\begin{align}
P_0[Z_t N_{\infty }] 
\le 
\liminf_{u \to \infty } 
P_0[Z_t \bar{N}_u] 
\le 
P_0[Z_t \bar{N}_t] , 
\end{align}
we see that $ P_0[N_{\infty }]<\infty $ and that 
$ (\xi_t:=\bar{N}_t-P_0[N_{\infty }|\cF_t]) $ 
is still a non-negative $ (\cF_t,P_0) $-supermartingale. 
Now the proof is completed by (i) and by Lemma \ref{lem: MtF/1+hXt}. 
\end{proof}

Finally, we proceed to prove Theorem \ref{thm: char of pos mart}. 

\begin{proof}[Proof of Theorem \ref{thm: char of pos mart}]
The necessity is immediate from Lemma \ref{lem: MtF/1+hXt}. 
Let us prove the sufficiency. 

Let $ (N_t) $ be a non-negative $ (\cF_t,P_0) $-martingale. 
Then, by Proposition \ref{prop: supermart decomp}, 
we have the decomposition \eqref{eq: supermart decomp}. 
Letting $ t=0 $, we have 
\begin{align}
N_0 = M_0(F) + P_0[N_{\infty }] + \xi_0 . 
\end{align}
Since $ N_0=M_0(F)=\sP[F] $ by the assumption, we have 
\begin{align}
P_0[N_{\infty }] = \xi_0 = 0 . 
\end{align}
This proves that $ N_t=M_t(F) $ $ P_0 $-almost surely, 
which completes the proof. 
\end{proof}


\end{document}